\documentclass[12pt]{article}

\usepackage{amsmath, amsthm, amsfonts, amssymb}
\usepackage{comment}
\usepackage[all]{xy}
\usepackage{graphics}
\usepackage{color}
\usepackage{fullpage}

\theoremstyle{plain}
\newtheorem{thm}{Theorem}[section]
\newtheorem{lem}[thm]{Lemma}
\newtheorem{prop}[thm]{Proposition}
\newtheorem{cor}[thm]{Corollary}

\theoremstyle{definition}
\newtheorem*{defn}{Definition}

\theoremstyle{remark}

\newcommand{\R}{\operatorname{Reach}}
\newcommand{\wt}{\operatorname{wt}}
\newcommand{\ceiling}[1]{\lceil{#1}\rceil}

\newcommand\peg[3]{\xymatrix{#1 \ar@/^2mm/[r]_{#2}& #3}}
\newcommand\peb{f}
\newcommand\twopeb{f_2}
\newcommand\dve{d_{\mathrm{ve}}}

\allowdisplaybreaks[1]

\title{Graph pegging numbers}
\author{
{\small
\begin{tabular}{cc}
{\large Geir Helleloid} & {\large Madeeha Khalid}\\[3pt]
Department of Mathematics & Department of Mathematics \\
The University of Texas at Austin & Institute of
Technology, Tralee \\
1 University Station C1200 & South Campus, Clash \\
Austin, TX 78712-0257 & Tralee, Co. Kerry, Ireland \\
\texttt{geir@math.utexas.edu} &
\texttt{madeeha.khalid@staff.ittralee.ie} \\[16pt]
{\large David Petrie Moulton}
       & {\large Philip Matchett Wood} \\[3pt]
Center for Communications Research & Department of Mathematics\\
805 Bunn Drive & Hill Center-Busch Campus\\
Princeton, NJ 08540 & Rutgers, The State University of New Jersey\\
\texttt{moulton@idaccr.org} & 110 Frelinghuysen Rd \\
& Piscataway, NJ 08854 \\
& \texttt{matchett@math.rutgers.edu}
\end{tabular}
}
}

\begin{document}
\maketitle

\begin{abstract}
In graph pegging, we view each vertex of a graph as a hole into which a
peg can be placed, with checker-like ``pegging moves'' allowed.
Motivated by well-studied questions in graph pebbling, we introduce two
pegging quantities.  The pegging number (respectively, the optimal
pegging number) of a graph is the minimum number of pegs such that for
every (respectively, some) distribution of that many pegs on the graph,
any vertex can be reached by a sequence of pegging moves. We prove
several basic properties of pegging and analyze the pegging number and
optimal pegging number of several classes of graphs, including paths,
cycles, products with complete graphs, hypercubes, and graphs of small
diameter.
\end{abstract}

\section{Introduction}\label{s:intro}

The mathematics of peg jumping originated in the games of Peg Solitaire
and Conway's Soldiers (see Berlekamp, Conway, and Guy~\cite{bcg}). Much
has been written about answering the classical questions of Peg
Solitaire and Conway's Soldiers in more general settings (see, for
example, Eriksen, Eriksson, and Eriksson~\cite{eee}).  In this paper,
however, we consider peg jumping -- or, as we call it, pegging -- on
graphs, and our questions are inspired by work in the theory of graph
pebbling.

Given a graph, we view each vertex as a hole into which one peg can be
placed.  A \emph{pegging move} consists of removing two pegs from
adjacent holes and placing one peg in a third, empty hole adjacent to
one of the first two holes.  In essence, one peg is jumping the other
and landing in the third hole (with the jumped peg being removed).  If
there are pegs in some of the vertices of the graph, we say that we can
\emph{peg to a vertex} if we can move a peg to that vertex with a
(possibly empty) sequence of pegging moves.

Conway's Soldiers can be recast in this general setting.  In that game,
the graph in question is embedded in the Cartesian plane with vertex set
$\mathbb{Z}^2$, and there is an edge between two vertices if their
Euclidean distance is 1.  Pegs are placed at all vertices in the lower
half-plane, and the challenge is to move a peg as far as possible into
the upper half-plane by peg-jumping.  Classically, however, a peg is
only allowed to jump in a straight line (horizontally or vertically)
over another peg, which is more restrictive than our pegging moves.  In
fact, whether or not all pegging moves are allowed, the best possible
solution moves a peg up four units.  This was proved in the classical
case by Conway using a weight argument; our results on graphs use an
extension of that weight argument and incidentally show that the optimum
solution to Conway's Soldiers cannot be improved upon using pegging
moves.

As mentioned, the pegging questions we consider are motivated by graph
pebbling (see Hurlbert~\cite{hurlbert} for a survey of graph pebbling).
The \emph{pegging number} of a graph is the minimum number of pegs so
that no matter how those pegs are distributed on the graph, we can peg
to any vertex.  In contrast, the \emph{optimal pegging number} of the
graph is the minimum number of pegs so that there is \emph{some} way to
distribute those pegs on the graph so that we can peg to any vertex.
These definitions mirror those of the pebbling number and the optimal
pebbling number of a graph.  In fact, we will use results on pebbling
numbers to prove results on pegging numbers; basic pebbling definitions
will be given when needed.

A formal definition of pegging and fundamental pegging lemmas appear in
Section~\ref{s:basic}.  In Section~\ref{s:classes}, we study the pegging
numbers of several classes of graphs, including paths and cycles.  In
Section~\ref{s:prod}, we move on to the Cartesian product of an
arbitrary graph $G$ with a complete graph $K_n$.  In particular,
Theorem~\ref{t:crossKn} relates the pegging number of $G \times K_n$ to
the pebbling number of $G$.  In Section~\ref{s:hypercube}, we apply this
result to compute the pegging number of the hypercube, and we also
obtain upper and lower bounds for the optimal pegging number of the
hypercube; this uses the theory of binary linear codes. Finally, in
Section~\ref{s:diameter}, we give an upper bound for the optimal pegging
number of a graph of diameter $2$, classify graphs with pegging number
at most 3, and give an upper bound for the pegging number of a graph of
diameter at most $3$.

The concept of pegging on a graph was introduced by the third author at
the 1994 University of Minnesota Duluth Research Experience for
Undergraduates (REU).  He and several of the students at the program
spent a weekend exploring basic properties of pegging.  Their results
are included here, along with our later work.  Further work on pegging
has been done by Wood~\cite{wood} and Levavi~\cite{lev}.  It should also
be mentioned that Niculescu and Niculescu~\cite{nn} independently
proposed the idea of pegging on graphs; however, their only result in
that direction is the first conclusion of Lemma~\ref{l:mon weight}.

\section{The Basics of Pegging}\label{s:basic}

A \emph{distribution} $D$ of pegs on a graph $G$ is any subset of
$V(G)$.  Through Lemma~\ref{l:mon weight}, $G$ may be any graph; after
Lemma~\ref{l:mon weight}, we only consider finite, simple graphs.  If
$u$ and $v$ are distinct, adjacent vertices in $D$, and $w$ is a vertex
adjacent to $v$ that is not in $D$, then the \emph{pegging move} $m =
\peg{u}{v}{w}$ replaces the distribution $D$ with the distribution $m(D)
= D \setminus \{u,v\} \cup \{w\}$.  Sometimes, to emphasize that the
conditions on $u$, $v$, and $w$ are met, we call $m$ a \emph{valid}
pegging move (on $D$).  The set $\{u,v\}$ is the \emph{source} of $m$,
and the vertex $w$ is the \emph{destination} of $m$.  If $M$ is a
sequence of pegging moves starting at $D$, then we write $M(D)$ for the
final distribution.  A vertex $t$ is \emph{reachable} from a
distribution $D$ if there is a finite sequence of pegging moves $M$ with
$t \in M(D)$.  The \emph{reach} of a distribution $D$, denoted $\R(D)$,
is the set of all vertices reachable from $D$.

The first of the four tools given in this section for analyzing pegging
is a weight argument adapted from the solution to Conway's Soldiers
given in~\cite{bcg}.  It is used to show that a given vertex is not in
the reach of a distribution.  Given a distribution $D$ on a graph $G$
and a vertex $t$, define the \emph{weight} of $D$ with respect to $t$ to
be 
\[ 
\wt_t(D) = \sum_{u \in D}{\sigma^{d(u,t)}}, 
\] 
where $\sigma = (\sqrt{5}-1)/2$ is the positive root of $x^2+x=1$, and
$d(u,t)$ is the distance in $G$ from $u$ to $t$.  Note that if $D$ or
$G$ is infinite, $\wt_t(D)$ may be infinite.

\begin{lem}[Monotonicity of Weight]\label{l:mon weight}
Let $D$ be a distribution on a graph $G$, and let $D'$ be a distribution
obtained from $D$ by a finite sequence of pegging moves.  Then
\[
\wt_t(D')\le \wt_t(D)
\]
for all $t \in G$.  If $\wt_t(D) < 1$, then $t \notin \R(D)$.
\end{lem}

\begin{proof}
Without loss of generality, we may assume that $D' = m(D)$, where $m =
\peg{u}{v}{w}$.  For any vertex $t$, we know $d(u,t) \le d(w,t) + 2$ and
$d(v,t) \le d(w,t)+1$.  Thus
\begin{eqnarray*}
\wt_t(D') &=& \wt_t(D) - \sigma^{d(u,t)} 
                   - \sigma^{d(v,t)} + \sigma^{d(w,t)} \\
&\le& \wt_t(D) - \sigma^{d(w,t)+2} 
             - \sigma^{d(w,t)+1} + \sigma^{d(w,t)} \\
&=& \wt_t(D) - (\sigma^2+\sigma-1)\sigma^{d(w,t)} \\
&=& \wt_t(D).
\end{eqnarray*}
The second claim follows from the observation if $t \in \R(D)$, then $t$
is in some distribution $D'$ obtained from $D$ by a finite sequence of
pegging moves.  But then $1 \le \wt_t(D') \le \wt_t(D)$.
\end{proof}

Combining Lemma~\ref{l:mon weight} with a short computation shows that
the optimum solution to Conway's Soldiers is still four units when all
pegging moves are allowed.  For the rest of the paper, we will only
consider finite, simple graphs $G$.  Our principal goal is to study the
following two pegging invariants.

\begin{defn}
The \emph{pegging number} of a graph $G$ is the smallest positive
integer $d$ such that \emph{every} distribution of size $d$ on $G$ has
reach $V(G)$.  The \emph{optimal pegging number} $p(G)$ of $G$ is the
smallest positive integer $d$ such that \emph{some} distribution of size
$d$ of $G$ has reach $V(G)$.  
\end{defn}

We can make some general observations about $P(G)$ and $p(G)$.
Obviously $p(G) \le P(G)$.  Let $|G|$ denote the \emph{order} of $G$,
that is, the number of vertices of $G$. If $G$ is disconnected, then
\[
P(G) = |G| - \min_C{(|C| - P(C))}
\]
and
\[
p(G) = \sum_C{p(C)},
\]
where the minimum and sum are taken over all connected components $C$ of
$G$.  Thus, we will focus primarily on connected graphs.

Let $\alpha(G)$ denote the \emph{independence number} of a graph $G$,
that is, the maximum cardinality of a set of pairwise non-adjacent
vertices of $G$.  Clearly $\alpha(G) \le P(G) \le |G|$.  Furthermore, if
$G$ has at least one edge, then $\alpha(G) + 1 \le P(G)$.  (In the
pebbling literature, graphs achieving equality in the corresponding
lower bound for the pebbling number are said to be ``of class~$0$'' or
``demonic'', so one might say that a graph $G$ with $P(G) = \alpha(G) +
1$ is ``of pegging class~$0$'' or ``devilish''.)  If no connected
component of $G$ is isomorphic to a star graph $K_{1,k}$, then $P(G) \le
|G| - 1$.  If $G$ has at least two vertices, then $p(G) \ge 2$, with
equality if and only if there are two adjacent vertices that dominate
$G$.  Finally, a simple application of Lemma~\ref{l:mon weight} proves
the following proposition, which was discovered by various
mathematicians at the 1994 University of Minnesota Duluth REU.

\begin{prop}\label{p:diameter}
If a graph $G$ has diameter $d$, then $P(G) \ge d$.
\end{prop}

\begin{proof}
Choose vertices $u_0$ and $u_d$ with $d(u_0, u_d) = d$, and let $u_0,
u_1, \dots, u_d$ be a path of length $d$ from $u_0$ to $u_d$.  Let $D =
\{u_2, \dots, u_d\}$.  Then
\[
\wt_{u_0}(D) = \sigma^2 + \sigma^3 + \cdots + \sigma^d < 1.
\]
Thus, by Lemma~\ref{l:mon weight}, $v_0 \notin \R(D)$ and $P(G) \ge d$.
\end{proof}

Our goal is to be able to compute the pegging numbers and optimal
pegging numbers of a variety of graphs.  Our first computational tool is
the weight argument given in Lemma~\ref{l:mon weight}.  Our remaining
three tools show that allowing the removal of pegs, ``stacking moves'',
and ``pebbling moves'' does not increase the reach of a distribution.
This helps us get upper bounds on pegging and optimal pegging numbers
because we can use all of these moves to show that the reach of a
distribution really is all of $V(G)$.  The definition of these moves and
the proof of the main result require a more sophisticated view of
pegging; in particular, we have to number the pegs so that they are
distinguishable and do not get ``mixed up'' when stacked on a vertex.

Given a graph $G$, let $S_G = \{u_i \; : \; u \in V(G), i \in
\mathbb{Z} \}$.  We interpret $u_i$ to indicate that peg $i$ is on
vertex $u$. Define a \emph{multi-distribution} $D$ of pegs on $G$ to be
a finite subset of $S_G$ with the property that if $u_i, v_i \in D$,
then $u = v$ (that is, peg $i$ can only be on one vertex at a time).

If $u$, $v$, and $w$ are distinct vertices, $v$ is adjacent to $u$ and
$w$, and $u_i$ and $v_j$ are in $D$, then the \emph{stacking move} $m =
\peg{u_i}{v_j}{w_i}$ sends the multi-distribution $D$ to the
multi-distribution $m(D) = D \setminus \{u_i, v_j\} \cup \{w_i\}$.  If
$u$ and $w$ are distinct, adjacent vertices, $i \neq j$, and $u_i$ and
$u_j$ are in $D$, then the \emph{pebbling move} $m =
\peg{u_i}{u_j}{w_i}$ sends the multi-distribution $D$ to the
multi-distribution $m(D) = D \setminus \{u_i, u_j\} \cup \{w_i\}$.  (See
Section~\ref{s:prod} for a discussion of pebbling.)  Finally, if $u_i$
is in $D$, then the \emph{removal move} $m$ (with respect to $u_i$)
sends the multi-distribution $D$ to the multi-distribution $m(D) = D
\setminus \{u_i\}$. 

We view distributions (in which the pegs happen to be labeled) as
multi-distributions in the obvious way, and we view pegging moves as
stacking moves in the obvious way.  To emphasize that a
multi-distribution is, in fact, a distribution, we may refer to it as a
\emph{proper} distribution.  If $D$ is a proper distribution, let
$\R_a(D)$ denote the set of vertices reachable from $D$ via all moves
(stacking, pebbling, and removal).  Note that $\R_a(D) \supseteq
\R(D)$.

Given a multi-distribution $D$, a \emph{move forest} of $D$ is a labeled
binary forest (that is, a disjoint union of labeled binary trees) with
the following three properties:
\begin{enumerate}
\item The label on each node is an element of $S_G$; multiple nodes may
have the same label.
\item The label on each leaf node is an element of $D$; no two leaf
nodes may have the same label.
\item Each interior node has a left child and a right child.  If an
interior node is labeled $w_i$, then either:
\begin{itemize}
\item the left and right children have labels of the form $u_i$ and
$v_j$, respectively, where $u$, $v$, and $w$ are distinct vertices, $v$
is adjacent to $u$ and $w$, and $i \neq j$, or
\item the left and right children have labels of the form $u_i$ and
$u_j$, respectively, where $u$ and $w$ are distinct, adjacent vertices
and $i \neq j$.
\end{itemize}
\end{enumerate}

A \emph{traversal} of a move forest is an ordering of the interior nodes
so that each interior node precedes its ancestors.  Each traversal $N =
(n_1, n_2, \dots)$ of a move forest corresponds to a sequence of valid
stacking and pebbling moves $M = (m_1, m_2, \dots)$ on $D$, where $m_r =
\peg{u_i}{v_j}{w_i}$ if node $n_r$, its left child, and its right child
are labeled $w_i$, $u_i$, and $v_j$ respectively, and $m_r =
\peg{u_i}{u_j}{w_i}$ if node $n_r$, its left child, and its right child
are labeled $w_i$, $u_i$, and $u_j$ respectively.  The following
proposition is clear.

\begin{prop}
Given a multi-distribution $D$, the above correspondence gives a
bijection between sequences of valid stacking and pebbling moves on $D$
and traversals of move forests on $D$.  
\end{prop}

\begin{thm} \label{t:reach}
Let $D$ be a proper distribution on a graph $G$.  Then $\R_a(D) = \R(D)$.
\end{thm}

\begin{proof}
The proof is by induction on the size of $D$.  The result is clear when 
$|D| = 1$, so assume it is true for all distributions of size less than
$d$, and let $D$ be a distribution of size $d$.  Let $t \in \R_a(D)
\setminus D$ and let $M$ be a sequence of stacking, pebbling, and
removal moves that puts a peg on $t$.  If we were to remove all the
removal moves from $M$, we would get a new sequence of stacking and
pebbling moves on $D$ that is valid and puts a peg on $t$ (the only
possible change is that some of the pegging moves in $M$ change to
non-pegging stacking moves).  So we may assume $M$ contains only
stacking and pebbling moves.  Let $F$ be the move forest corresponding
to $M$; note that $M$ has an interior node labeled $t_a$ for some $a \in
\mathbb{Z}$.  By the inductive hypothesis, it suffices to show that
there exists a valid pegging move $m'$ on $D$ so that $t \in
\R_a(m'(D))$.

Next we show that we can assume $M$ has no pebbling moves.  Suppose $m =
\peg{u_i}{u_j}{w_i}$ is a pebbling move in $M$, and let $n$ be the node
in $F$ corresponding to $m$.  Since $D$ is proper, the left child
of $n$ (labeled $u_i$) or the right child of $n$ (labeled $u_j$) is not
a leaf node.  We may assume without loss of generality that the left
child $n'$ of $n$ is not a leaf node.  (If the right child of $n$ is not
a leaf node, we can swap the left and right subtrees of $n$ and change
every occurrence of $i$ in the labels of $n$ and its ancestors to a $j$
to form a new move forest that puts a peg on $t$ in which the left child
of $n$ is not a leaf node.)  If the right child of $n'$ is labeled $v_k$
for some vertex $v \neq w$, replace the left subtree of $n$ by the right
subtree of $n'$ and change every occurrence of $i$ in the labels of $n$
and its ancestors to a $k$.  If the right child of $n'$ is labeled
$w_k$, replace $n$ and its subtree by the right subtree of $n'$ and
change every occurrence of $i$ in the labels of the (former) ancestors
of $n$ to a $k$.  In either case, we have a new move forest that puts a
peg on $t$ and has fewer pebbling moves than $M$.  We can repeat this
operation until all pebbling moves have been removed.  Therefore, we may
assume that $M$ has no pebbling moves.

Finally we handle stacking moves.  Let $m$ be the first move in $M$
whose target vertex is not in $D$.  Write $m = \peg{u_i}{v_j}{w_i}$.
Note that by our choice of $m$, the distribution $D$ does contain $u_k$
and $v_l$ for some $k$ and $l$, and $m' = \peg{u_k}{v_l}{w_k}$ is a
valid pegging move on $D$.  We only need to show that $t \in
\R_a(m'(D))$.  In fact, we will construct a new move forest $F'$ with a
node labeled $t_b$ for some $b \in \mathbb{Z}$ and an interior node
labeled $w_k$ whose left and right children, respectively, are leaf
nodes labeled $u_k$ and $v_l$. Then any traversal of $F'$ starting with
this interior node puts peg $b$ at vertex $t$ and shows $t \in
\R_a(m'(D))$.

Let $n$ be the node in $F$ corresponding to $m$.  Let $A$ be the left
subtree of $n$ and let $B$ be the right subtree of $n$.  Let $C$ be the
set of leaf nodes in $F$ that are not leaves of $A$ or $B$.  Detach $A$
and $B$ from $n$, replacing $A$ by a single leaf node labeled $u_k$ and
replacing $B$ by a single leaf node labeled $v_l$.

There are some situations in which we will reattach $A$ or $B$ to a different
part of the tree, and we describe those situations in this paragraph.
If $u_k$ is not the label of a leaf node in $B$, or $v_j$ is not the
label of a leaf node in $A$, we may assume without loss of generality
that $u_k$ is not the label of a leaf node in $B$.  If $u_k$ is the
label of a leaf node in $C$, replace that leaf node by $A$.  If $v_l$ is
the label of a leaf node in $C$ or $A$, replace that leaf node by $B$.

Finally, choose any traversal of the current forest.  Visiting each
interior node in order, change the subscript on the label of the node to
match the subscript on the label of its left child (if they are
different).  Let the resulting forest be $F'$; we claim it is a move
forest on $D$.

Property 1 in the definition of a move forest is clearly satisfied.
Property 2 must be satisfied since we only added leaf nodes labeled
$u_k$ and $v_l$ and removed any other leaf nodes with those labels.
Property 3 was also enforced in the previous paragraph.  Furthermore,
$F'$ has the interior node $n$ labeled $w_k$ whose left and right
children, respectively, are leaf nodes labeled $u_k$ and $v_l$.
Finally, the node labeled $t_a$ in $F$ must exist in $F'$, but possibly
with the label changed to $t_b$.  As mentioned, any traversal of $F'$
starting with $n$ shows $t \in \R_a(m'(D))$, proving the theorem.
\end{proof}

\begin{cor}[Monotonicity of Reach] \label{c:mon reach}
Let $D' \subset D$ be two distributions on a graph $G$.  Then $\R(D')
\subseteq \R(D)$.  If $D$ is a distribution of size $d$ and $\R(D)
\neq V(G)$, then $P(G) > d$.  If $\R(D) \neq V(G)$ for every
distribution $D$ of size $d$, then $p(G) > d$.
\end{cor}

\section{Paths, Cycles, and Joins}\label{s:classes}

In this section, our goal is to compute the pegging and optimal pegging
numbers of several simple classes of graphs.  We will use $K_n$ to
denote the complete graph on $n$ vertices; $P_n$ to denote the path on
$n$ vertices; and $C_n$ to denote the cycle on $n$ vertices.  For
convenience, we will label the vertices of both $P_n$ and $C_n$ by $v_1,
v_2, \dots, v_n$ with $v_i$ adjacent to $v_{i+1}$ for $1 \le i \le n-1$.
The complete graphs are simple to analyze: $P(K_n) = p(K_n) = 2$ for $n
\ge 2$.  As for cycles, $P(C_3) = 2$ and $P(C_4) = 3$, so we turn to
cycles on five or more vertices.

\begin{thm}\label{t:cycle}
For $n \ge 5$, the pegging number of the cycle $C_n$ is $P(C_n) = n-2$.
\end{thm}

\begin{proof}
In every distribution of $n-2$ pegs on $C_n$, for each of the two holes,
there are two adjacent pegs with one of them adjacent to the hole, so
that each hole is the destination for some move.  Hence $P(C_n) \le
n-2$.

On the other hand, let $D = \{v_2, v_3, \dots, v_{n-2}\}$.  Suppose
$v_n \in \R(D)$, and let $M$ be a minimum-length sequence of pegging
moves with $v_n \in M(D)$.  By minimality, the last move of $M$ is the
first move placing a peg at $v_n$, and by symmetry, we may assume
that the last move of $M$ is $\peg{v_{n-2}}{v_{n-1}}{v_n}$.  Then $M$ is
a valid pegging sequence on the distribution $D' = \{v_2, v_3, \dots,
v_{n-2}\}$ on the path $P_n$.  This says $v_n \in \R(D')$ while
$\wt_{v_n}(D') < 1$, contradicting Lemma~\ref{l:mon weight}.  So $\R(D)
\neq V(C_n)$ and $P(C_n) > n-3$.
\end{proof}

\begin{thm}\label{t:opt cycle}
For $n \ge 3$, the optimal pegging number of the cycle $C_n$ is $p(C_n)
= \ceiling{n/2}$. 
\end{thm}

\begin{proof}
Let $D$ be the distribution $\{v_2, v_3, v_5, v_7, v_9, \dots, v_n\}$ if
$n$ is odd and $\{v_2, v_3, v_5, v_7, v_9, \dots, v_{n-1}\}$ if $n$ is
even.  Then $R(D) = V(C_n)$ and $|D| = \ceiling{n/2}$, so $p(C_n) \le
\ceiling{n/2}$.

Suppose there is a distribution $D$ of $\ceiling{n/2} - 1$ pegs with
$\R(D) = V(C_n)$.  $D$ consists of blocks of consecutive pegs
alternating with blocks of consecutive empty vertices.  If $D$ has a
block of three or more empty vertices, then (by symmetry) we may assume
that $D$ is contained in the distribution $\{v_2, v_3, \dots,
v_{n-2}\}$, whose reach is not $V(C_n)$ by the proof of 
Theorem~\ref{t:cycle}.  By Corollary~\ref{c:mon reach}, this contradicts
our choice of $D$, so $D$ does not contain a block with three or more
empty vertices.

$D$ has strictly more empty vertices than pegs.  Empty vertices appear
in blocks of one or two, so there are more blocks of two empty vertices
than blocks of two or more pegs.  In particular, since $D$ must have at
least one block of two or more pegs, there are at least two blocks of
two empty vertices.  It follows that there exist empty vertices 
$u$, $v$, $w$, and $x$ so that $u$ and $v$ are adjacent, $w$ and $x$ are
adjacent, and the pegs in one of the two components $G_1$ and $G_2$ of
$C_n \setminus \{(u,v), (w,x)\}$ only come in blocks of one.  Without
loss of generality, let it be $G_1$, with $v$ and $w$ adjacent to $G_1$.
No moves among the pegs in $G_1$ are possible, since none are adjacent,
until a sequence of moves on the pegs in $D' = D \cap G_2$ puts a peg on
$v$ or $w$.  However, $\wt_v(D') < 1$ and $\wt_w(D') < 1$, so $v,w \notin
\R(D')$.  Thus $v$ and $w$ are not in $\R(D)$, contradicting our choice
of $D$.  The result follows.
\end{proof}

Using the optimal pegging number of the cycle $C_n$ and the following
proposition, we can compute $P(P_n)$ and $p(P_n)$.

\begin{prop}\label{p:spanning}
If $H$ is a spanning subgraph of $G$, then $p(G) \le p(H)$ and $P(G) \le
P(H)$. 
\end{prop}

\begin{proof}
Any move made on a distribution on $H$ can be made on the same
distribution on $G$, so the reach of a distribution on $G$ contains the
reach of the same distribution on $H$.  The desired inequalities follow.
\end{proof}

For $1 \le n \le 3$, we clearly have $P(P_n) = n$, so we now consider
the pegging numbers of paths of order at least $4$.

\begin{thm}\label{t:path}
For $n \ge 4$, the pegging number of the path $P_n$ is $P(P_n) = n-1$.
\end{thm}

\begin{proof}
It is clear that $\R(D) = V(P_n)$ for any distribution $D$ of $n-1$ pegs
on $P_n$.  On the other hand, the diameter of $P_n$ is $n-1$, so by
Proposition~\ref{p:diameter}, we have $P(P_n) \ge n-1$.  Thus $P(P_n) =
n-1$.
\end{proof}

Obviously $p(P_n) = n$ for $n = 1$ and $2$, so we now consider larger
values of $n$.

\begin{thm}\label{t:opt path}
For $n \ge 3$, the optimal pegging number of the path $P_n$ is $p(P_n) =
\ceiling{n/2}$.  
\end{thm}

\begin{proof}
Since $P_n$ is a spanning subgraph of $C_n$,
Proposition~\ref{p:spanning} shows that $p(P_n) \ge p(C_n) =
\ceiling{n/2}$.  On the other hand, let $D$ be the distribution $\{v_2,
v_3, v_5, v_7, \dots, v_n\}$ if $n$ is odd and $\{v_2, v_3, v_5, v_7,
\dots, v_{n-1}\}$ if $n$ is even.  Then $\R(D) = P_n$, so $p(P_n) \le
\ceiling{n/2}$.  Hence equality holds.  
\end{proof}

We close this section by calculating the pegging and optimal pegging
numbers of joins.  The \emph{join} of two graphs $G$ and $H$, denoted
$G + H$, has vertex set $V(G) \cup V(H)$ and edge set $E(G) \cup E(H)
\cup \{(g,h) : g \in G, h \in H\}$.  Recall that $\alpha(G)$ is the
independence number of a graph $G$.

\begin{thm}\label{t:join}
Given any two graphs $G$ and $H$ each having at least one vertex, the
join $G+H$ satisfies $p(G+H) = 2$ and 
\[ 
P(G+H) = \alpha(G+H) + 1 = \max(\alpha(G), \alpha(H))+1.
\] 
\end{thm}

\begin{proof}
Clearly $2$ is a lower bound for $p(G+H)$.  To achieve this lower bound,
place one peg on a vertex of $G$ and one peg on a vertex of $H$.  It
is easy to see that the reach of this distribution is all of $G+H$,
giving $p(G+H) = 2$.

Take $a = \max(\alpha(G), \alpha(H))$.  Note that $\alpha(G+H) = a$, so
$P(G+H) \ge a + 1$.  To show that $P(G+H)$ is exactly $a+1$, first
observe that if a distribution has one peg on a vertex of $G$ and one
peg on a vertex of $H$, then the reach is the entire graph.  Let $D$
be a distribution of $a+1$ pegs on $G+H$ in which either all the pegs
are on $G$ or all the pegs are on $H$.  Without loss of generality,
let all the pegs be on $G$.

If $a = 1$, then $G$ and $H$ are both complete graphs and so $G+H$ is
complete.  Since $|D| = 2$, it follows that $\R(D) = G+H$.  If $a > 1$,
then $|D| \ge 3$.  Since $|D| > \alpha(G)$, there are adjacent vertices
$u$ and $v$ in $D$.  Choose any vertex $w$ in $H$ and let $m =
\peg{u}{v}{w}$. Then in the distribution $m(D)$ there is a peg on a
vertex of $H$ and there remains a peg on a vertex of $G$, so every
vertex is in the reach of $m(D)$.  Hence we have $\R(D) = G+H$ for any
distribution of $a+1$ pegs on $G+H$, and we obtain $P(G+H) =
\max(\alpha(G), \alpha(H)) + 1$.  
\end{proof}

\begin{cor}
For any complete multi-partite graph $G$, we have $p(G) = 2$, and $P(G)$
is one more than the size of the largest partite set of $G$. 
\end{cor}

\begin{proof}
This follows immediately from Theorem~\ref{t:join} and the fact that $G$
is the join of one partite set and either another partite set or a
complete multi-partite subgraph. 
\end{proof}

\section{Products with Complete Graphs}\label{s:prod}
The Cartesian product of two graphs $G$ and $H$, denoted $G \times H$,
has
vertex set
\[
V(G \times H) = \{(g,h) : g \in V(G), h \in V(H)\},
\]
and there is an edge between vertices $(g, h)$ and $(g', h')$ if $g =
g'$ and $h$ and $h'$ are adjacent in $H$ or if $h = h'$ and $g$ and $g'$
are adjacent in $G$.  We often view a Cartesian product as consisting of
$|H|$ copies of $G$ with additional edges between different copies.
Computing the pebbling numbers of a Cartesian product is one of the
primary directions in the pebbling literature (see \cite{hurlbert}).  In
this section, we study the pegging numbers of the Cartesian product $G
\times K_n$.  The main result of this section, Theorem~\ref{t:crossKn},
provides an upper bound for the pegging number based on pebbling
numbers.  (A formal definition of pebbling is given below.)  Before
presenting Theorem~\ref{t:crossKn}, we first study the pegging and
optimal pegging numbers of the Cartesian product $K_m \times K_n$.

\begin{prop}
For positive integers $m$ and $n$, not both 1, the pegging number of
$K_m \times K_n$ is 
\[ 
P(K_m \times K_n) = \alpha(K_m \times K_n) + 1 = \min\{m,n\} + 1.
\] 
\end{prop}

\begin{proof}
Without loss of generality, we may assume $m \ge n$.  Since $\alpha(K_m
\times K_n) = n$, we know $P(K_m \times K_n) \ge n+1$.  Let $D$ be any
distribution of $n+1$ pegs on $K_m \times K_n$ and let $t$ be any target
vertex; we may assume $t \notin D$.  

Let $A$ be the copy of $K_m$ containing $t$.  If $A$ has pegs on two
vertices $u$ and $v$, then $\peg{u}{v}{t}$ puts a peg on $t$.  If $A$
only has one peg, say on $u$, then some other copy $B$ of $K_m$ has pegs
on two vertices $v$ and $w$.  Let $x$ be the vertex in $A$ adjacent to
$w$.  Then the moves $\peg{v}{w}{x}$ (which may be a stacking move) and
(if $x$ is not $t$) $\peg{x}{u}{t}$ (which may be a pebbling move) put a
peg on $t$.

Finally, if $A$ has no pegs, either two copies of $K_m$ each have (at
least) two pegs or one copy $B$ has pegs at three vertices $u$, $v$, and
$w$.  In the first case, one move puts a peg on $A$, reducing to the
case when $A$ has one peg.  In the second case, let $x$ be the vertex in
$B$ adjacent to $t$.  If $x$ equals $u$, $v$, or $w$, one move puts a
peg on $t$.  If not, then $\peg{u}{v}{x}$ and $\peg{w}{x}{t}$ put a peg
on $t$.  Thus $t \in \R(D)$ and $P(K_m \times K_n) = n+1$.  
\end{proof}

\begin{prop}\label{p:opt product}
For positive integers $m$ and $n$, the optimal pegging number of $K_m
\times K_n$ is 
\[
p(K_m \times K_n) = \left\{ \begin{array}{r@{\quad\quad}l}
1 & \text{if } mn = 1, \\
2 & \text{if } mn > 1 \text{ and } \min\{m,n\} \le 2, \\
3 & \text{otherwise.}
\end{array} \right.
\]
\end{prop}

\begin{proof}
Without loss of generality, we may assume $m \ge n$.  The result is
trivial for $n = 1$.  When $n = 2$, place pegs on two vertices in one
copy of $K_n$.  Any vertex can be reaching in one move, showing that
$p(K_m \times K_n) = 2$.

Finally, suppose $n \ge 3$.  Any placement of two pegs on $K_m \times
K_n$ leaves a copy of $K_m$ and a copy of $K_n$ with no pegs, so we
cannot peg to their common vertex.  Thus $p(K_m \times K_n) > 2$.  On
the other hand, place pegs at three vertices $u$, $v$, and $w$ in one
copy $A$ of $K_m$ and let $t$ be any vertex.  If $t$ is in $C$ or just
adjacent to $u$, $v$, or $w$, one move puts a peg at $t$.  Otherwise,
let $x$ be the vertex in $C$ adjacent to $t$.  Then $\peg{u}{v}{x}$ and
$\peg{w}{x}{t}$ put a peg on $t$.  So $p(K_m \times K_n) = 3$.
\end{proof}

We now move on to the main result of this section.  Consider a graph $G$
with some number of pebbles placed on each vertex.  A \emph{pebbling
move} consists of removing two pebbles from one vertex and placing one
pebble on an adjacent vertex.  The \emph{pebbling number} $\peb(G)$
(respectively, the \emph{two-pebbling number} $\twopeb(G)$) of a graph
$G$ is the smallest positive integer $p$ so that, given any distribution
of $p$ pebbles on $G$ and any vertex $t$, there exists a sequence of
pebbling moves that places a pebble at $t$ (respectively, two pebbles at
$t$).  Note that $\peb(G) \ge |G|$ and that $\peb(G)$ and $\twopeb(G)$
are taken to be infinite if $G$ is disconnected.

\begin{thm}\label{t:crossKn}
For any graph $G$ and positive integer $n$, the pegging number of $G
\times K_n$ satisfies $P(G \times K_n) \le \twopeb(G)$.  Furthermore,
$P(G \times K_2) \le \max\{\peb(G), |G|+1\}$. 
\end{thm}

\begin{proof}
The result is trivial when $G$ is disconnected, so we may assume $G$ is
connected.  If $v$ is a vertex in $G$, let $K_v$ denote the copy of
$K_n$ in $G \times K_n$ consisting of $n$ copies of $v$.  If $w$ is a
vertex in $K_n$, let $G_w$ denote the copy of $G$ in $G \times K_n$
consisting of $|G|$ copies of $w$.  Let $\pi : V(G \times K_n) \rightarrow
V(G)$ be the projection map.

Let $D$ be any distribution of $\twopeb(G)$ pegs on $G \times K_n$ and
let $t$ be any target vertex.  Let $g = \pi(t)$ and view the
multi-distribution $\pi(D)$ as a distribution of pebbles on $G$.  There
are $\twopeb(G)$ pebbles in $\pi(D)$, so there is a sequence of pebbling
moves on $G$ that places two pebbles on $g$.  For each pebbling move on
$G$, perform a corresponding stacking or pebbling move on $G \times K_n$
(that is, given a pebbling move from $h \in G$ to $k \in G$, perform a
stacking move from two vertices in $K_h$ to a vertex in $K_k$ or a
pebbling move from a vertex in $K_h$ to a vertex in $K_k$).  This puts
two pegs on $K_g$.  One more stacking or pebbling move puts a peg on
$t$.  By Theorem~\ref{t:reach}, $t \in \R(D)$.  Thus $P(G \times K_n) \le
\twopeb(G)$.

Now let $p = \max\{\peb(G), |G|+1\}$, let $D$ be any distribution of
$p$ pegs on $G \times K_2$, and let $t$ be any target vertex (not in
$D$).  Let $g = \pi(t)$ and view the multi-distribution $\pi(D)$ as a
distribution of pebbles on $G$.  Let $u$ be the other vertex in $K_g$
besides $t$.

If $K_g$ does not contain any pegs, then we can move a pebble to $g$,
since there are $p \ge \peb(G)$ pebbles on $G$.  Perform corresponding
stacking or pebbling moves on $G \times K_2$ to move a peg to $K_g$,
always taking the target of each move to be in the copy of $G$
containing $t$.  (Since we only stack pegs in this copy, we will never
have to make a pebbling move using two pegs in the other copy of $G$.
The final move will place a peg at $t$.

On the other hand, suppose $K_g$ does have a peg (on $u$).  Since $p >
|G|$, some vertex $h \neq g$ has at least two pebbles.  Choose the
nearest such $h$ to $g$ and let $g = x_0, x_1, x_2, \dots, x_k = h$
be a minimum length path $P$ from $g$ to $h$.  Note that $x_1, \dots,
x_{k-1}$ each have at most one pebble.  If they each have exactly one
pebble, perform the sequence of pebbling moves $(m_k, m_{k-1}, \dots,
m_1)$, where $m_i$ removes two pebbles from $x_i$ and puts one pebble on
$x_{i-1}$.  Any corresponding sequence of stacking and pebbling moves on
$G \times K_2$ puts another peg in $K_g$, and $t$ can be reached with
(at most) one more stacking or pebbling move.  

If, however, not all of the vertices $x_1, \dots, x_{k-1}$ have a
pebble, let $j \ge 1$ be the least index such that $x_j$ does not have a
pebble.  Then we can move a pebble to $x_j$.  Note that $x_0, x_1,
\dots, x_{j-1}$ each have only one pebble, so if any of these pebbling
moves involved one of these vertices, the first such move would have to
land a pebble on $x_i$ for some $i$ with $0 \le i < j$.  After this
move, $x_1, \dots, x_{i-1}$ have one pebble and $x_i$ has two pebbles,
and as in the previous paragraph, we can reach $t$. Hence we may assume
that none of $x_0, x_1, \dots, x_{j-1}$ are involved in moving a pebble
to $x_j$.  As before, using an appropriate corresponding sequence of
stacking and pebbling moves, we can get a peg to the vertex $u$ of
$K_{x_j}$ that is adjacent to the vertex $v$ of $K_{x_{j-1}}$ having a
peg. Now we can make a pegging move from $u$ over $v$ to a vertex of
$K_{x_{j-2}}$ (or to $t$ for $j=1$), and again we are reduced to the
case of the previous paragraph, from which we can reach $t$.

In every case, we can reach $t$ by making stacking and pebbling moves
from the distribution $D$, so Theorem~\ref{t:reach} shows that $P(G
\times K_2) \le p$, proving the theorem.  
\end{proof}

\section{Hypercubes}\label{s:hypercube}

For each positive integer $n$, the hypercube $Q_n$ is the graph $K_2
\times K_2 \times \cdots \times K_2$, a product of $n$ copies of $K_2$.
For convenience, we label the vertices of each factor $K_2$ with $0$ and
$1$, and we label each of the $2^n$ vertices of $Q_n$ with the
corresponding binary $n$-tuple. An early result of Fan Chung in pebbling
literature is that the pebbling number of $Q_n$ is $2^n$
(see~\cite{chung}).  Using this fact and Theorem~\ref{t:crossKn}, we can
compute the pegging numbers of hypercubes.

\begin{thm}
For any positive integer $n$, we have $P(Q_n) = 2^{n-1} + 1$.
\end{thm}

\begin{proof}
We have $\peb(Q_{n-1}) = 2^{n-1} = |Q_{n-1}|$ (see \cite{chung}), so
Theorem~\ref{t:crossKn} gives $P(Q_n) = P(Q_{n-1} \times K_2) \le
2^{n-1} + 1$.  On the other hand, the set consisting of the $2^{n-1}$
vertices of $Q_n$ for which the sum of the $n$ coordinates is even is an
independent set, so we have $P(Q_n) \ge \alpha(Q_n) + 1 \ge 2^{n-1} +
1$.  This proves the theorem. 
\end{proof}

It seems harder to compute the exact optimal pegging number of a
hypercube, but we can give upper and lower bounds that are correct up to
a polynomial factor.

\begin{lem}\label{l:opt_peg_prod}
Given any distributions $D$ and $E$ on graphs $G$ and $H$, respectively,
\[
\R(D) \times \R(E) \subseteq \R(D \times E),
\]
where $D \times E$ is viewed as a distribution on $G \times H$.
\end{lem}

\begin{proof}
Let $t \in \R(D)$ and $u \in \R(E)$.  We must show that $(t,u) \in \R(D
\times E)$.  For each $v \in E$, make the moves on the subgraph $G
\times \{v\}$ using $D \times \{v\}$ needed to put a peg on $(t,v)$.
Now there is a peg on each vertex in $\{t\} \times E$, so we can make
moves in $\{t\} \times H$ to put a peg on $(t,u)$.
\end{proof}

\begin{cor}\label{c:opt_peg_prod}
For any graphs $G$ and $H$ we have $p(G \times H) \le p(G) p(H)$.
\end{cor}

Now we obtain our bounds on the optimal pegging numbers of cubes. The
results and methods are similar to those of Moews in the case of
pebbling (see~\cite{moews}).  The upper bound was originally obtained by
Lenhard L. Ng in 1994.

\begin{thm}\label{t:opt_peg_cube}
For any positive integer $n$, 
\[
\left(\sqrt{5}-1\right)^n \le p(Q_n) \le (2n)^{3/2}
\left(\sqrt{5}-1\right)^n.
\]
\end{thm}

\begin{proof}
To prove the lower bound, let $D$ be any distribution on $Q_n$ with
$\R(D) = V(Q_n)$. By Lemma~\ref{l:mon weight}, $\wt_t(D) \ge 1$
for all $t \in V(Q_n)$.  Summing over all the vertices of $Q_n$ shows
\begin{align*}
2^n & \le \sum_{t \in Q_n} \wt_t(D)\\
&= \sum_{t \in Q_n} \sum_{v \in D} \wt_t(v)\\
&= \sum_{v \in D} \sum_{t \in Q_n} \wt_t(v).
\end{align*}
The contribution of the peg at $v \in D$ to this sum is
\begin{align*}
\sum_{t \in Q_n} \wt_t(v) &= \sum_{i = 0}^n \mathop{\sum_{t \in
Q_n}}_{d(v,t) =
i} \wt_t(v)\\
&= \sum_{i = 0}^n{\binom{n}{i} \sigma^i}\\
&= (1+\sigma)^n\\
&= \phi^n,
\end{align*}
with $\phi = (1+\sqrt{5})/2$.  Inserting this into the previous inequality
shows 
\[
|D| \phi^n \ge 2^n
\]
and 
\[
|D| \ge \left( \frac{2}{\phi} \right)^n = \left(\sqrt{5}-1\right)^n,
\]
which gives the lower bound.

In order to obtain the upper bound, we fix an integer $n \ge 2$, let $r
= \left\lceil n/(\phi+2)\right\rceil$, and let $m = n - r$. Then there
is a binary linear code $C$ in $Q_m$ with covering radius at most $r$
and dimension at most $d = m(1 - H(r/m)) + (3/2) \lg m + 1$, with $H(x)
= -x \lg x - (1-x) \lg (1-x)$ and $H(0) = H(1) = 0$ (see~\cite{cohen};
``$\lg$'' here means logarithm to the base $2$). That is, there is a
(linear) subset $C$ of $Q_m$ with at most $2^d$ vertices such that every
vertex of $Q_m$ has distance at most $r$ from some vertex of $C$. 

Now let $D$ be the distribution $Q_r \times C$ on $Q_r \times Q_m = Q_n$
and take any vertex $(u,v) \in Q_r \times Q_m$.  We have a vertex $w \in
C$ at distance at most $r$ from $v$ in $Q_m$.  This means that some copy
of $Q_r$ in $Q_m$ contains both $v$ and $w$, so that the
sub-distribution $Q_r \times \{w\} \subseteq D$ and the vertex $(u,v)$
both lie in some copy of $Q_r \times Q_r$.  Now note that 
\[
\R(Q_r \times Q_0) = Q_r \times Q_r; 
\]
this follows immediately by induction on $r \ge 1$, using
Lemma~\ref{l:opt_peg_prod}. Thus $(u,v)$ is in the reach of the
sub-distribution $Q_r \times \{w\}$, which means that all of $Q_n$ is in
the reach of the distribution $D$.

This gives us $p(Q_n) \le |Q_r \times C| \le 2^{r+d}$, and we have
\begin{align*}
r + d & = r + m\left(1 - H\left(\frac{r}{m}\right)\right)
            + \frac{3}{2} \lg m + 1\\
      & = \frac{3}{2} \lg m + n + 1 - (n-r)H\left(\frac{r}{n-r}\right).
\left(\frac{m-r}{m}\right)^{m-r}.
\end{align*}
Now $H(x)$ is increasing on the interval
$[0\mathinner{\ldotp\ldotp}1/2]$, so we have
\begin{align*}
r + d & \le \frac{3}{2} \lg m + n + 1
            - (n-r)H\left(\frac{n/(\phi+2)}{n-n/(\phi+2)}\right)\\
      &  =  \frac{3}{2} \lg m + n + 1
            - (n-r)H\left(\frac{1}{\phi+1}\right)\\
      & \le \frac{3}{2} \lg m + n + 1
            - \left(n-\frac{n}{\phi + 2}-1\right)H\left(\phi^{-2}\right)\\
      & \le \frac{3}{2} \lg n + n + 1
         + \left[\lg\left(1-\frac{1}{\phi+2}\right) + H(\phi^{-2})\right]
         - n\left(\frac{\phi + 1}{\phi + 2}\right)H\left(\phi^{-2}\right)\\
      & \le \frac{3}{2} \lg n + n + 1 + \frac{1}{2}
            - n\left(\frac{\phi^2}{\phi + 2}\right)
            (2\phi^{-2}\lg\phi + \phi^{-1}\lg\phi)\\
      &  =  \frac{3}{2} \lg(2n) + n - n\lg\phi\\
      &  =  \frac{3}{2} \lg(2n) + n\lg\left(\sqrt{5}-1\right). \end{align*}
We conclude that 
\[
p(Q_n) \le (2n)^{3/2} \left(\sqrt{5} - 1\right)^n, 
\]
as desired. 
\end{proof}

\section{Graphs of Small Diameter}\label{s:diameter}

We can prove sharp bounds on the pegging numbers and optimal pegging
numbers of graphs of small diameter.  Our results will require merely
that the graphs have vertex-edge diameter (defined below) at most some
value $d$, which is weaker than requiring the diameter to be at most
$d$.

We define the distance between a vertex and an edge of a graph to be the
minimum distance between the vertex and an endpoint of the edge.  We
denote the diameter of a graph $G$ by $d(G)$ and the \emph{vertex-edge
diameter} (the maximum distance between a vertex and an edge) of a
non-null graph $G$ by $\dve(G)$.  (Recall that a graph is non-null if at
has at least one edge).  Note that a non-null graph $G$ of diameter $d$
has $d-1 \le \dve(G) \le d$.  Notice also that the vertex-edge diameter
of a graph, unlike the diameter, can be increased by adding an edge; for
instance adding an edge to $K_{1,3}$ increases the vertex-edge diameter
from $1$ to $2$.

Let $G$ be a non-null graph with $\dve(G) \le 1$.  From the definition
of optimal pegging number, we have $p(G) = 2$, since if we put pegs on
any two adjacent vertices, we can jump to any other vertex in one move.
We claim that we also have $P(G) = \alpha(G) + 1$. For if we are given
any target vertex $t \in G$ and a distribution $D$ of $\alpha(G) + 1$
vertices not containing $t$, then $D$ contains two adjacent vertices,
one of which is adjacent to $t$, and we may jump the peg on one of these
vertices over the other one to $t$.

Now let $G$ be a graph of order at least $2$ and radius $1$, so that it
has a vertex $u$ adjacent to every other vertex.  Then $G$ has optimal
pegging number $2$, since $u$ and any other vertex together have reach
equal to $G$. We claim that we also have $P(G) = \alpha(G) + 1$.  If $G$
is complete, the result is trivial, so we may assume $\alpha(G) \ge 2$.
Given any distribution $D$ of at least $\alpha(G) + 1$ vertices, $D$
contains two adjacent vertices, so we can jump one over the other to
reach $u$ (if $u$ is not already in $D$).  Since we have $|D| \ge 3$, we
can still jump from another vertex over $u$ to any desired target (not
already in $D$). This proves the claim.

We now obtain sharp bounds on the optimal pegging number and pegging
number of graphs of vertex-edge diameter $2$.

\begin{prop}\label{p:ld2}
The optimal pegging number of any graph $G$ with $\dve(G) = 2$
is at most~$4$. \end{prop}

\begin{proof}
Let $G$ be a graph with $\dve(G) = 2$.  Now $G$ must have two disjoint
edges $(u, v)$ and $(w, x)$, and we let $D$ be $\{u, v, w, x\}$. Take
any target vertex $t \in G$.  If $t$ is in $D$ or is adjacent to a
vertex of $D$, then we can reach $t$ in at most one move.  Otherwise,
either $u$ or $v$, say $u$, has a common neighbor $y$ with $t$.  We may
assume that $y$ is not adjacent to $w$ or $x$ and then, similarly, $w$,
say, must have a common neighbor $z$ with $y$.  Now we can jump from $v$
over $u$ to $y$, from $x$ over $w$ to $z$, and from $z$ over $y$ to $t$,
reaching the target.  Thus $G$ is in the reach of $D$ and we have $p(G)
\le 4$. 
\end{proof}

\begin{figure}
\begin{center}
\scalebox{.75}{\includegraphics{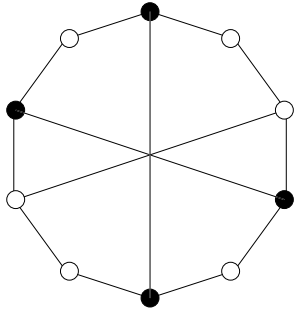}}
\caption{The graph above has vertex-edge diameter $2$. The distribution
consisting of the black vertices has reach equal to the entire graph, so
the optimal pegging number is at most $4$. It can be shown that no
distribution of $3$ vertices can reach every vertex, so the optimal
pegging number is exactly $4$.} \label{dpm0}
\end{center}
\end{figure}

The authors experimented by computing the optimal pegging numbers of a
number of small graphs of vertex-edge diameter $2$.  It seemed
emperically that such graphs that do not contain a pair of adjacent
vertices dominating the graph have optimal pegging number $3$, however
one example with optimal pegging number $4$ is given in
Figure~\ref{dpm0}.  While there do not seem to be any simple examples of
graphs with diameter $2$ and optimal pegging number $4$, there is the
Hoffman--Singleton graph.  This is the unique (up to isomorphism)
$7$-regular graph of order $50$, diameter $2$, and girth $5$ (so any two
non-adjacent vertices have exactly one common neighbor)
(see~\cite{hoffsing}).

Suppose we have some distribution of three vertices $u, v, w$ on the
Hoffman--Singleton graph.  There are two possibilities up to symmetry:
either $u$ and $w$ are both adjacent to $v$, or $u$, $v$, $x$, $w$, and
$y$ are the five vertices of a pentagon listed in order, for some
vertices $x$ and $y$. In the first case, it is not possible to make two
moves in sequence, so the total reach of the distribution consists of
$20$ vertices:\ the $3$ vertices $u, v, w$, plus $6$ neighbors each of
$u$ and $w$ that are not $v$, plus $5$ neighbors of $v$ that are not $w$
or $u$. In the second case, it is possible to make at most two moves in
sequence, but only if the first move is $\peg{u}{v}{x}$ or
$\peg{v}{u}{y}$. Thus, the reach of the second distribution consists of
the $5$ vertices $\{u,v,x,w,y\}$, together with their $25$ other
neighbors, making the total number of vertices in the reach $30$. Hence
the optimal pegging number of the Hoffman--Singleton graph is $4$.

For the pegging numbers of graphs of diameter $2$, we obtain results
similar to those of Clarke, Hochberg, and Hurlbert (see~\cite{chh}).

\begin{figure}
\begin{center}
%\scalebox{.75}{%\input{dpm_fig1.pstex_t}
\begin{picture}(0,0)%
\includegraphics{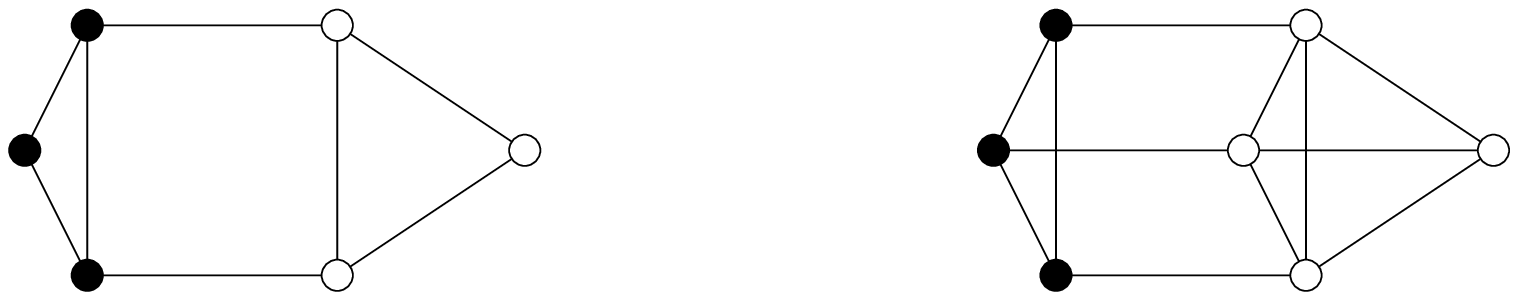}%
\end{picture}%
\setlength{\unitlength}{3947sp}%
\begingroup\makeatletter\ifx\SetFigFont\undefined%
\gdef\SetFigFont#1#2#3#4#5{%
  \reset@font\fontsize{#1}{#2pt}%
  \fontfamily{#3}\fontseries{#4}\fontshape{#5}%
  \selectfont}%
\fi\endgroup%
\begin{picture}(7583,2469)(151,-2224)
\put(7651,-961){\makebox(0,0)[lb]{\smash{\SetFigFont{12}{14.4}{\familydefault}{\mddefault}{\updefault}{\color[rgb]{0,0,0}$t$}%
}}}
\put(5176,
14){\makebox(0,0)[lb]{\smash{\SetFigFont{12}{14.4}{\familydefault}{\mddefault}{\updefault}{\color[rgb]{0,0,0}$v_1$}%
}}}
\put(6676,
89){\makebox(0,0)[lb]{\smash{\SetFigFont{12}{14.4}{\familydefault}{\mddefault}{\updefault}{\color[rgb]{0,0,0}$u_1$}%
}}}
\put(1201,-2161){\makebox(0,0)[lb]{\smash{\SetFigFont{14}{16.8}{\familydefault}{\mddefault}{\updefault}{\color[rgb]{0,0,0}(a)}%
}}}
\put(6001,-2161){\makebox(0,0)[lb]{\smash{\SetFigFont{14}{16.8}{\familydefault}{\mddefault}{\updefault}{\color[rgb]{0,0,0}(b)}%
}}}
\put(3001,-961){\makebox(0,0)[lb]{\smash{\SetFigFont{12}{14.4}{\familydefault}{\mddefault}{\updefault}{\color[rgb]{0,0,0}$t$}%
}}}
\put(526,
14){\makebox(0,0)[lb]{\smash{\SetFigFont{12}{14.4}{\familydefault}{\mddefault}{\updefault}{\color[rgb]{0,0,0}$v_1$}%
}}}
\put(451,-1486){\makebox(0,0)[lb]{\smash{\SetFigFont{12}{14.4}{\familydefault}{\mddefault}{\updefault}{\color[rgb]{0,0,0}$v_2$}%
}}}
\put(2026,
89){\makebox(0,0)[lb]{\smash{\SetFigFont{12}{14.4}{\familydefault}{\mddefault}{\updefault}{\color[rgb]{0,0,0}$u_1$}%
}}}
\put(151,-586){\makebox(0,0)[lb]{\smash{\SetFigFont{12}{14.4}{\familydefault}{\mddefault}{\updefault}{\color[rgb]{0,0,0}$v_3$}%
}}}
\put(1651,-1561){\makebox(0,0)[lb]{\smash{\SetFigFont{12}{14.4}{\familydefault}{\mddefault}{\updefault}{\color[rgb]{0,0,0}\quad
$u_2$}%
}}}
\put(4801,-586){\makebox(0,0)[lb]{\smash{\SetFigFont{12}{14.4}{\familydefault}{\mddefault}{\updefault}{\color[rgb]{0,0,0}$v_3$}%
}}}
\put(5101,-1486){\makebox(0,0)[lb]{\smash{\SetFigFont{12}{14.4}{\familydefault}{\mddefault}{\updefault}{\color[rgb]{0,0,0}$v_2$}%
}}}
\put(6076,-586){\makebox(0,0)[lb]{\smash{\SetFigFont{12}{14.4}{\familydefault}{\mddefault}{\updefault}{\color[rgb]{0,0,0}$u_3$}%
}}}
\put(6451,-1561){\makebox(0,0)[lb]{\smash{\SetFigFont{12}{14.4}{\familydefault}{\mddefault}{\updefault}{\color[rgb]{0,0,0}$u_2$}%
}}}
\end{picture}
%}
\caption{Let $G$ be a graph with $\alpha(G) = \alpha \ge 2$ and with
pegging number $\alpha +2$.  If $G$ has vertex-edge diameter $2$, then
$G$ contains (a) as an induced subgraph, and if $G$ has diameter $2$,
then $G$ contains (b) as an induced subgraph. A vertex is black or white
according to whether it lies in the distribution in the proof of
Theorem~\ref{t:ld2}.} \label{dpm1} 
\end{center} 
\end{figure}

\begin{thm}\label{t:ld2}
For any integer $\alpha \ge 2$, the maximum pegging number of graphs of
vertex-edge diameter $2$ and independence number $\alpha$ is $\alpha +
2$, and this bound can be achieved by a graph of diameter $2$.
Furthermore, given a distribution of at least this size for such a
graph, we always can reach any target with at most $3$ moves.  Finally,
any such graph with pegging number exactly $\alpha + 2$ contains the
graph of Figure~\ref{dpm1}(a) (respectively, Figure~\ref{dpm1}(b)) as an
induced subgraph. 
\end{thm}

\begin{figure}
\begin{center}
\begin{picture}(0,0)%
\includegraphics{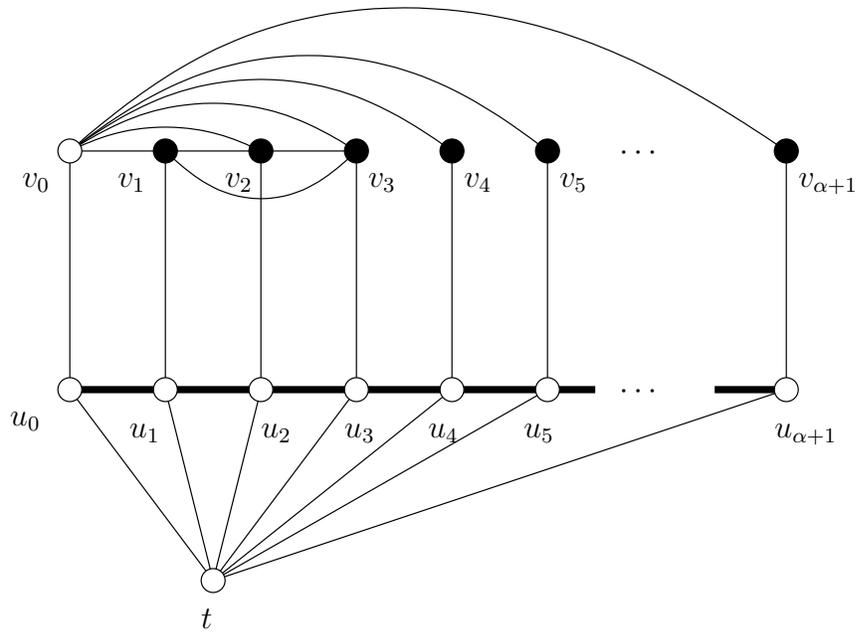}%
\end{picture}%
\setlength{\unitlength}{3947sp}%
\begingroup\makeatletter\ifx\SetFigFont\undefined%
\gdef\SetFigFont#1#2#3#4#5{%
  \reset@font\fontsize{#1}{#2pt}%
  \fontfamily{#3}\fontseries{#4}\fontshape{#5}%
  \selectfont}%
\fi\endgroup%
\begin{picture}(4958,3961)(226,-3410)
\put(4051,-1915){\makebox(0,0)[lb]{\smash{\SetFigFont{12}{14.4}{\familydefault}{\mddefault}{\updefault}{\color[rgb]{0,0,0}$\cdots$}%
}}}
\put(4051,-415){\makebox(0,0)[lb]{\smash{\SetFigFont{12}{14.4}{\familydefault}{\mddefault}{\updefault}{\color[rgb]{0,0,0}$\cdots$}%
}}}
\put(1426,-3361){\makebox(0,0)[lb]{\smash{\SetFigFont{12}{14.4}{\familydefault}{\mddefault}{\updefault}{\color[rgb]{0,0,0}$t$}%
}}}
\put(226,-2086){\makebox(0,0)[lb]{\smash{\SetFigFont{12}{14.4}{\familydefault}{\mddefault}{\updefault}{\color[rgb]{0,0,0}$u_0$}%
}}}
\put(976,-2161){\makebox(0,0)[lb]{\smash{\SetFigFont{12}{14.4}{\familydefault}{\mddefault}{\updefault}{\color[rgb]{0,0,0}$u_1$}%
}}}
\put(1801,-2161){\makebox(0,0)[lb]{\smash{\SetFigFont{12}{14.4}{\familydefault}{\mddefault}{\updefault}{\color[rgb]{0,0,0}$u_2$}%
}}}
\put(2326,-2161){\makebox(0,0)[lb]{\smash{\SetFigFont{12}{14.4}{\familydefault}{\mddefault}{\updefault}{\color[rgb]{0,0,0}$u_3$}%
}}}
\put(2851,-2161){\makebox(0,0)[lb]{\smash{\SetFigFont{12}{14.4}{\familydefault}{\mddefault}{\updefault}{\color[rgb]{0,0,0}$u_4$}%
}}}
\put(3451,-2161){\makebox(0,0)[lb]{\smash{\SetFigFont{12}{14.4}{\familydefault}{\mddefault}{\updefault}{\color[rgb]{0,0,0}$u_5$}%
}}}
\put(5026,-2161){\makebox(0,0)[lb]{\smash{\SetFigFont{12}{14.4}{\familydefault}{\mddefault}{\updefault}{\color[rgb]{0,0,0}$u_{\alpha+1}$}%
}}}
\put(301,-586){\makebox(0,0)[lb]{\smash{\SetFigFont{12}{14.4}{\familydefault}{\mddefault}{\updefault}{\color[rgb]{0,0,0}$v_0$}%
}}}
\put(901,-586){\makebox(0,0)[lb]{\smash{\SetFigFont{12}{14.4}{\familydefault}{\mddefault}{\updefault}{\color[rgb]{0,0,0}$v_1$}%
}}}
\put(1576,-586){\makebox(0,0)[lb]{\smash{\SetFigFont{12}{14.4}{\familydefault}{\mddefault}{\updefault}{\color[rgb]{0,0,0}$v_2$}%
}}}
\put(2476,-586){\makebox(0,0)[lb]{\smash{\SetFigFont{12}{14.4}{\familydefault}{\mddefault}{\updefault}{\color[rgb]{0,0,0}$v_3$}%
}}}
\put(3076,-586){\makebox(0,0)[lb]{\smash{\SetFigFont{12}{14.4}{\familydefault}{\mddefault}{\updefault}{\color[rgb]{0,0,0}$v_4$}%
}}}
\put(3676,-586){\makebox(0,0)[lb]{\smash{\SetFigFont{12}{14.4}{\familydefault}{\mddefault}{\updefault}{\color[rgb]{0,0,0}$v_5$}%
}}}
\put(5176,-586){\makebox(0,0)[lb]{\smash{\SetFigFont{12}{14.4}{\familydefault}{\mddefault}{\updefault}{\color[rgb]{0,0,0}$v_{\alpha+1}$}%
}}}
\end{picture}

\caption{
The graph above both has diameter and vertex-edge diameter equal to $2$
and has independence number $\alpha$. The vertices $u_i$ form a complete
graph on $\alpha + 2$ vertices (the thick lines denote the many edges
joining these vertices), and the $\alpha +1$ vertices that are colored
black represent a distribution of pegs from which it is not possible to
reach $t$. In fact, the only pegging moves that leave adjacent pegs
involve moving a peg to $v_0$, and then no more useful moves can be
made. Thus the pegging number of this graph is at least $\alpha + 2$. }
\label{dpm2} 
\end{center} 
\end{figure}

\begin{proof}
First note that the graph in Figure~\ref{dpm2} has both diameter and
vertex-edge diameter equal to $2$ and has independence number $\alpha$,
but has pegging number at least $\alpha + 2$.

Now suppose we have a graph $G$ with $\dve(G) = 2$, a target vertex $t$
of $G$, and a distribution $D$ of at least $\alpha(G) + 2$ pegs on $G$.
Then there are two adjacent vertices $u$ and $v$ of $D$, and one of
them, say $u$, has distance at most $2$ from $t$. We may assume the
distance is $2$, so $t$ and $u$ have a common neighbor $w$. Now move a
peg from $v$ over $u$ to $w$. The remaining distribution still has two
adjacent pegs (and we may assume that neither is $w$), so we may, as
before, make a second move to get a peg adjacent to $w$.  Finally, with
a third move, we jump this peg over $w$ to $t$.  Thus, given a graph $G$
of vertex-edge diameter $2$ and any distribution of at least $\alpha(G)
+ 2$ pegs, any target can be reached with at most $3$ moves, and we have
$P(G) \le \alpha(G) + 2$.

Finally, suppose we have a graph $G$ with $\dve(G) = 2$ and $P(G) =
\alpha(G) + 2$.  Then we have a target vertex $t$ of $G$ and a
distribution $D$ of $\alpha(G) + 1$ pegs whose reach does not contain
$t$. Notice that no vertex of $D$ can be equal to or adjacent to $t$, or
we would be able to get a peg to $t$ in at most two moves, as in the
previous paragraph.  Similarly, there cannot exist two disjoint pairs of
adjacent vertices of $D$, or we would be able to get to $t$ in at most
three moves. Let $v_1$ and $v_2$ be two adjacent vertices of $D$; then
$\left(D \setminus \{v_1\}\right) \cup \{t\}$ has more than $\alpha(G)$
vertices and so contains two adjacent vertices, which must lie in $D
\setminus \{v_1\}$.  One of these two vertices must be $v_2$, lest $D$
contain two disjoint adjacent pairs, so we have $v_2$ adjacent to some
vertex $v_3 \in D \setminus \{v_1\}$. Similarly, $\left(D \setminus
\{v_2\}\right) \cup \{t\}$ contains an adjacent pair of vertices, and
this pair must consist of $v_1$ and $v_3$. Thus $v_1, v_2, v_3$ are all
adjacent, and these are the only adjacencies in $D$.

Now one of $v_1$ and $v_2$, say $v_1$, must have a common neighbor $u_1$
with $t$, and $u_1$ cannot be adjacent to any vertex $v \in D \setminus
\{v_1\}$, for then we would be able to jump from $v_2$ (or $v_3$ if $v$
is $v_2$) over $v_1$ to $u_1$ and then from $v$ over $u_1$ to $t$.
Similarly, one of $v_2$ and $v_3$, say $v_2$, has a common neighbor
$u_2$ with $t$, and $u_2$ cannot be adjacent to any vertex of $D
\setminus \{v_2\}$. Finally, the set $\left(D \setminus \{v_1,
v_2\}\right) \cup \{u_1, u_2\}$ contains two adjacent vertices, and
these must be $u_1$ and $u_2$, so $u_1$ and $u_2$ are adjacent.  Thus
the subgraph of $G$ induced by $\{t, u_1, u_2, v_1, v_2, v_3\}$ is
isomorphic to the graph of Figure~\ref{dpm1}(a). In addition, if $G$ has
diameter $2$, then there must, similarly, be a vertex $u_3$ adjacent to
$t, v_3, u_1, u_2$ and not adjacent to $v_1, v_2$, so $G$ contains an
induced subgraph isomorphic to the graph of Figure~\ref{dpm1}(b).
\end{proof}

We now consider graphs with pegging number at most $3$.  First of all,
$K_0$ and $K_1$ are the only graphs with pegging numbers $0$ and $1$,
respectively.  Also, since the pegging number of a non-null graph is
greater than its independence number, we see that $2 K_1$ and $K_n$ for
$n \ge 2$ are the only graphs with pegging number $2$.  We can now use
Theorem~\ref{t:ld2} to classify graphs with pegging number $3$.  Recall
that we sometimes identify $D$ with the subgraph of $G$ spanned by $D$.

\begin{cor}\label{pegnum3}
A graph $G$ has pegging number $3$ if and only if $G$ is isomorphic to
$3 K_1$, or $G$ has independence number $2$ and the following condition
fails:
\begin{itemize}
\item[{\rm($\ast$)}] $G$ is spanned by $C \cup D$ with $C$ and $D$
complete subgraphs of $G$, and either $G$ is isomorphic to $2K_1$ or
$2K_2$, or $G$ has vertices $t \in C$ and $v_1, v_2, v_3 \in D$, such
that $t$ is not adjacent to any vertex of $D$, and every vertex of $C$
is adjacent to at most one of $v_1, v_2, v_3$. 
\end{itemize} 
\end{cor}

\begin{proof}
We know that $3K_1$ has pegging number $3$, while other graphs with
independence number at least $3$ have pegging number at least $4$. Also,
$2K_1$ and $2K_2$ have pegging numbers $2$ and $4$, respectively, while
all other graphs satisfying the forbidden condition ($\ast$) have
pegging number at least $4$, since $t$ is not in the reach of the
distribution $\{v_1, v_2, v_3\}$.  Thus, in order to prove the
corollary, we may assume that we have $\alpha(G) = 2$ and $P(G) > 3$,
and we must show that the forbidden condition ($\ast$) holds.

We first claim that $G$ is spanned by a disjoint union of complete
subgraphs.  If $G$ has diameter at least $3$, then take two vertices $x$
and $y$ at distance at least $3$ from each other. Then every vertex of
$G \setminus \{x, y\}$ is adjacent to exactly one of $x$ and $y$, and
$x$ together with its neighbors form a complete subgraph of $G$, as do
$y$ and its neighbors.  Thus we may assume that $G$ has diameter at most
$2$.  Now by (the proof of) Theorem~\ref{t:ld2}, we know that $G$ has
vertices $t, v_1, v_2, v_3$ such that $t$ is not in the reach of the
distribution $\{v_1, v_2, v_3\}$.  Thus no neighbor of $t$ is adjacent
to more than one of $v_1, v_2, v_3$, and any two neighbors of $t$ have a
common non-neighbor among $v_1, v_2, v_3$ and are, therefore, adjacent.
Then $t$ and its neighbors form a complete subgraph of $G$, as do the
non-neighbors of $t$.  This proves the claim.

Now $G$ is spanned by a disjoint union of complete subgraphs $C$ and
$D$, and there are vertices $t, v_1, v_2, v_3 \in G$ with $t$ not in the
reach of $\{v_1, v_2, v_3\}$.  We may assume that $t$ is in $C$, and
then at least two of $v_1, v_2, v_3$, say $v_1$ and $v_2$, are in $D$,
since otherwise we would be able to reach $t$ in one move.  First
suppose that $v_3$ is in $D$. If any vertex $x \in G$ were either in $C$
and adjacent to, say, $v_1$ and $v_2$, or in $D$ and adjacent to $t$,
then we could move $\peg{v_3}{v_1}{x}, \peg{v_2}{x}{t}$.  Thus neither
of these occurs, and ($\ast$) holds. Now suppose that $v_3$ is in $C$.
If any vertex $x \in G$ were either in $C$ and adjacent to, say, $v_1$,
or in $D$ and adjacent to $v_3$, then we could move $\peg{v_2}{v_1}{x},
\peg{x}{v_3}{t}$. Therefore, neither $v_1$ nor $v_2$ is adjacent to any
vertex of $C$, and $v_3$ is not adjacent to any vertex of $D$.  Now if
$D$ has a vertex $y$ besides $v_1$ and $v_2$, then ($\ast$) holds with
$v_3$ replaced by $y$ and $t$ replaced by $v_3$. Thus we may assume that
$D$ has order $2$, so there are no edges between the vertices of $C$ and
$D$.  Then either $G$ is isomorphic to $2K_2$, or $C$ has order at least
$3$, and in the latter case, interchanging the roles of $C$ and $D$
shows that ($\ast$) holds.  This proves the corollary. 
\end{proof}

Finally, we obtain a sharp upper bound on the pegging numbers of graphs
of vertex-edge diameter $3$.

\begin{thm}\label{t:ld3}
For any integer $\alpha \ge 2$, the maximum pegging number of graphs of
vertex-edge diameter $3$ and independence number $\alpha$ is $2\alpha +
1$, and this bound can be achieved by a graph of diameter $3$.
Furthermore, given a distribution of at least this size for such a
graph, we can always reach any target in at most $7$ moves. 
\end{thm}

First we need three lemmas that will be used several times in the proof
of this theorem.  The proofs of these results also introduce several
ideas that will appear in the proof of the theorem.

\begin{lem}\label{l:getp4}
Let $G$ be a connected graph with $|G| > 3$ and $|G| \ge 2\alpha(G) -
1$. Then $G$ contains a subgraph isomorphic to the path $P_4$. 
\end{lem}

\begin{proof}
Let $v$ be a vertex of $G$ of degree at least $2$, with neighbors $u$
and $w$.  If every neighbor of $v$ had degree $1$, then $G$ would be
isomorphic to the complete bipartite graph $K_{1,|G|-1}$.  But then we
would have 
\[ 
|G| \ge 2\alpha(G) - 1 = 2(|G|-1) - 1 = 2|G| - 3, 
\] 
which gives $|G| \le 3$, contrary to hypothesis.  Hence some neighbor of
$v$, say $w$, has degree at least $2$ and is adjacent to a vertex
besides $v$.  If $w$ is adjacent to $u$, then, as $G$ is connected and
contains some vertex besides $u, v, w$, one of these three vertices, say
$w$, is adjacent to a fourth vertex.  Thus, in any case, we may assume
that $w$ is adjacent to some vertex $x$ not equal to $u$ or $v$. Then
$G$ contains the path $uvwx$ of order $4$, proving the lemma.
\end{proof}

\begin{lem}\label{l:usep4}
Let $G$ be a graph with $\dve(G) = 3$, let $t \in G$ be a target vertex,
and let $D \subset G$ be a distribution not containing $t$, such that
$D$ contains a subgraph $P$ isomorphic to $P_4$.  Then we can reach a
neighbor of $t$ in at most $3$ moves using only the pegs of $P$.
\end{lem}

\begin{proof}
Label the vertices of $P$ in order as $u, v, w, x$.  Now either $v$ or
$w$, say $v$, has distance at most $3$ from $t$.  If this distance is
$1$, we are done, and if it is $2$, we can jump from $u$ over $v$ to a
neighbor of $t$, so we may assume that $v$ has distance $3$ from $t$.
Then $v$ has a neighbor $y$ adjacent to a neighbor $z$ of $t$, and we
can make the sequence of pegging moves $\peg{u}{v}{y}, \peg{x}{w}{v},
\peg{v}{y}{z}$, reaching a neighbor of $t$ in at most $3$ moves using
only the pegs of $P$, as desired. 
\end{proof}

\begin{lem}\label{l:2alpha}
Let $G$ be a graph with $\dve(G) = 3$, let $t \in G$ be a target vertex,
let $D \subset G$ be a distribution not containing $t$, and let $S
\subseteq G \setminus D$ be an independent set satisfying the following
conditions:
\begin{itemize}
\item No vertex of $S$ is adjacent to a vertex of $D$. \item No neighbor
of a vertex of $S$ both has distance $2$ from $t$ and is adjacent to a
vertex of a component of $D$ isomorphic to $K_3$. 
\item We have $|D| + d + 2|S| \ge 2\alpha(G) + 1$, with $d$ being the
number of isolated vertices of the subgraph $D$ of $G$. 
\end{itemize}
Then we can reach a neighbor of $t$ by making at most $3$ moves using the pegs
of $D$. 
\end{lem}

\begin{proof}
We may assume that every vertex of $D$ has distance at least $2$ from
$t$ and that every non-isolated vertex of $D$ has distance at least $3$
from $t$, since otherwise we can reach a neighbor of $t$ in at most one
move. For each component $C$ of $D$ that is isomorphic to $K_3$, take a
vertex $v_C$ of $C$ at distance $3$ from $t$ and a neighbor $u_C$ (in $G
\setminus D$) of $v_C$ at distance $2$ from $t$. For any component $C$
of $D$, we let $\hat{C}$ denote either the subgraph of $G$ induced by $C
\cup \{u_C\}$, for $C \cong K_3$, or $C$, otherwise.  We will refer to
$\hat{C}$ for $C \cong K_3$ as an ``extended subgraph''.

If some such vertex $u_C$ were adjacent to a vertex $w$ of $D$ besides
$v_C$, then we could jump from some vertex of $C$ over $v_C$ to $u_C$
and then from $w$ over $u_C$ to a neighbor of $t$, so we may assume that
this is not the case.  Similarly, if $u_C$ were adjacent to a vertex
$u_{C'}$ corresponding to another component $C'$, then we could jump
over $v_C$ to $u_C$, over $v_{C'}$ to $u_{C'}$, and from $u_{C'}$ over
$u_C$ to a neighbor of $t$, so we may assume that this does not occur,
either. Thus the extended subgraphs have independence number $2$, and
the components of $E = D \cup \{u_C \,|\, C \cong K_3 \text{ a component
of } D\}$ consist of the graphs $\hat{C}$, for $C$ a component of $D$.
Notice also that, by the hypotheses, no vertex of $S$ is adjacent to a
vertex of $E$.

We now have 
\begin{align*}
|D| &  >  2\alpha(G) - d - 2|S|\\
    & \ge 2\alpha(E \cup S) - d - 2|S|\\
    &  =  2\alpha(E) + 2|S| - d - 2|S|\\
    &  =  2\alpha(E) - d
\end{align*}
and, therefore, 
\[
\sum_{C \text{ a component of } D} (|C| - 2\alpha(\hat{C})) =
|D| - 2\alpha(E) > - d
\]
and 
\[
\mathop{\sum_{C \text{ a component of } D}}_{|C| > 1}
(|C| - 2\alpha(\hat{C})) > 0. 
\]
Hence we have a component $K$ of $D$ with $|K| > 1$ and $|K| >
2\alpha(\hat{K})$.  Now $K$ cannot be isomorphic to $K_3$, since then
$\hat{K}$ would be one of the extended subgraphs above and would have
independence number $2$.  Thus $K = \hat{K}$ has order at least $4$.
Now by Lemma~\ref{l:getp4}, $K$ contains a subgraph isomorphic to $P_4$,
and by Lemma~\ref{l:usep4}, we can get a peg to a neighbor of $t$ in at
most $3$ moves using only the pegs of $K \subseteq D$.  This proves the
lemma. 
\end{proof}

\begin{proof}[Proof of Theorem~\ref{t:ld3}.]
First, we exhibit a family of graphs achieving the maximum pegging
number.  Given a value $\alpha \ge 2$, we define a graph $G$ as follows:
The vertices of $G$ are $\{v_{ij} \,|\, 2 \le i \le \alpha \text{ and }
1 \le j \le 2, \text{ or } (i,j) \in \{(1,1), (\alpha, 3), (\alpha,
4)\}\} \cup \{u_{ij} \,|\, 1 \le i,j \le \alpha, i \ne j\}$.  All
vertices with first coordinate $i$ are adjacent, and, in addition,
$u_{ij}$ and $u_{ji}$ are adjacent.  Let $t$ be $v_{11}$ and $D$ be the
set of all other $v_{ij}$. Then $G$ is a graph of diameter $3$ and
independence number $\alpha$, and $D$ is a set of $2\alpha$ vertices
whose reach does not include $t$. Thus we have $P(G) \ge 2\alpha(G) +
1$.

Now suppose that we have a graph $G$ with $\alpha(G) = \alpha$ and
$\dve(G) = 3$, together with a target vertex $t$ and a distribution $D$
of at least $2\alpha + 1$ vertices. We wish to show that we can reach
$t$ from $D$ in at most $7$ moves. We may assume that $D$ does not
contain $t$.  Suppose that $D$ contains a neighbor $t'$ of $t$ or a
non-isolated (in $D$) vertex adjacent to a neighbor $t'$ of~$t$.  By
making one move if necessary, we may reduce to the first case with $|D|
\ge 2\alpha$.  Now if $t'$ had a neighbor in $D$, we could reach $t$ in
one (more) move, so assume not. Take $S = \{t'\}$ and $D' = D \setminus
S$.  Then we have 
\[ 
|D'| + 2|S| = |D| + 1 \ge 2\alpha + 1, 
\] 
so by Lemma~\ref{l:2alpha}, we know that we can make at most $3$ pegging
moves not involving $t'$ to get a peg to a neighbor of $t'$.  Then we
can jump this peg over $t'$ to $t$ to reach $t$ in a total of at most
$5$ moves.  Thus we may assume that every vertex of $D$ has distance at
least $2$ from $t$ and that every non-isolated vertex of $D$ has
distance at least $3$ from $t$.

As in the proof of Lemma~\ref{l:2alpha}, for each component $C$ of $D$
that is isomorphic to $K_3$, we take a vertex $v_C \in C$ at distance
$3$ from $t$ and a neighbor $u_C$ of $v_C$ at distance $2$ from $t$.  As
before, we denote by $\hat{C}$ the extended subgraph $C \cup \{u_C\}$
for such components $C$ (and $C$ itself for $C$ not isomorphic to
$K_3$). Suppose that some vertex $u_C$ is adjacent to a vertex $w \in D
\setminus C$. Then we can jump from a vertex of $C$ over $v_C$ to $u_C$
and from $w$ over $u_C$ to a neighbor $t'$ of $t$.  Now let $D'$ be the
current distribution minus $t'$, with $d'$ isolated vertices, and let
$S$ be $\{t'\}$.  If $t'$ has a neighbor in $D'$, then we can get to $t$
in one more move, so we may assume it does not.  We have $d' \ge 1$,
since the remaining vertex of $C$ is isolated in $D'$, and we obtain 
\[
|D'| + d' + 2|S| \ge |D| - 3 + 1 + 2 \ge 2\alpha + 1.  
\] 
Now we can apply Lemma~\ref{l:2alpha} to $t', D', S$ to show that we can
reach a neighbor of $t'$ with at most $3$ moves of the pegs of $D'$,
giving a total of at most $5$ moves.  Finally, we can jump from this
neighbor of $t'$ over $t'$ to $t$, reaching $t$ in at most $6$ moves.
Thus we may assume that no vertex $u_C$ is adjacent to any vertex of $D
\setminus C$.

Now suppose that some vertex $u_C$ is adjacent to a vertex $w$ of $C$
besides $v_C$.  As before, we can jump from the other vertex of $C$ over
$v_C$ to $u_C$ and then from $w$ over $u_C$ to a neighbor $t'$ of $t$.
Let $D'$ be the current distribution minus $t'$ and $S$ be $\{t', w\}$.
Then $w$ is not adjacent to any vertex of $D'$, and we may assume that
$t'$ is not either, or we would be able to get to $t$ in one more move.
Now if some neighbor $x$ of $w$ were adjacent to a vertex of some $C'$,
then we could peg to $x$, from $x$ over $w$ to $u_C$, and from $u_C$
over $t'$ to $t$, so we may assume that this is not the case.  Thus no
neighbor of $w$ violates the second condition of Lemma~\ref{l:2alpha},
and we have 
\[
|D'| + 2|S| \ge |D| - 3 + 4 \ge 2\alpha + 2,
\]
so we can apply Lemma~\ref{l:2alpha} as before to show that we can reach
$t$ in a total of at most $6$ moves.  Hence we may assume that no vertex
$u_C$ is adjacent to any vertex of $D$ besides $v_C$.

Finally, suppose two vertices $u_C$ and $u_{C'}$ with $C \ne C'$ are
adjacent.  Then we can jump over $v_C$ to $u_C$, over $v_{C'}$ to
$u_{C'}$, and from $u_{C'}$ over $u_C$ to a neighbor $t'$ of $t$.  Let
$D'$ be the current distribution minus $t'$, with $d'$ isolated
vertices, and let $S$ be $\{t'\}$.  As before, we have $d' \ge 2$, we
may assume that $t'$ has no neighbor in $D'$, and we calculate 
\[
|D'| + d' + 2|S| \ge |D| - 4 + 2 + 2 \ge 2\alpha + 1.
\]
By applying Lemma~\ref{l:2alpha}, we can reach $t$ in at most $4$ more
moves, giving a total of at most $7$ moves.  Thus we may assume that no
two vertices $u_C$ and $u_{C'}$ are adjacent.

Hence, as in the proof of Lemma~\ref{l:2alpha}, the extended subgraphs
all have independence number $2$, and the components of the graph $E = D
\cup \{u_C \,|\, C \cong K_3 \text{ a component of } D\}$ consist of the
graphs $\hat{C}$ for $C$ a component of $D$. Again, as in the proof of
Lemma~\ref{l:2alpha}, we have 
\begin{align*}
|D| & \ge 2\alpha(G) + 1\\
    & \ge 2\alpha(E) + 2\alpha(\{t\}) + 1\\
    &  =  2\alpha(E) + 3
\end{align*}
and, therefore, 
\begin{equation}\label{e:comps}
\sum_{C \text{ a component of } D} (|C| - 2\alpha(\hat{C})) \ge 3.
\end{equation}
Thus we have a component $K$ of $D$ with $|K| > 2\alpha(\hat{K})$.
Again, $K$ cannot be isomorphic to $K_3$, as $\hat{K}$ would then be an
extended subgraph with independence number~$2$. Therefore $K$ must equal
$\hat{K}$ and must have order at least $4$. By Lemma~\ref{l:getp4}, $K$
then contains a subgraph isomorphic to $P_4$.

For any (not necessarily induced) subgraph $P$ of $D$ isomorphic to
$P_4$, let $E'$ be obtained from $D \setminus P$ by extending each
component $C$ isomorphic to $K_3$ with a vertex $u_C$ as before.  For
each component $C$ of $D \setminus P$, we let $\check{C}$ be $C \cup
\{u_C\}$ if $C$ is isomorphic to $K_3$ and $C$ otherwise, paralleling
the notation $\hat{C}$ for components of $D$.  Suppose that some
extended subgraph $\check{C}$ does not have independence number $2$ or
is not a component of $E'$ or that some component $C$ of $D \setminus P$
satisfies $|C| > 3$ and $|C| \ge 2\alpha(C) - 1$.  Then either some
vertex $u_C$ is adjacent to a vertex of $D \setminus P$ besides the
corresponding $v_C$ or is adjacent to a vertex $u_{C'}$ for $C' \ne C$,
or, by Lemma~\ref{l:getp4}, $C$ must contain a subgraph isomorphic to
$P_4$. As above, or by applying Lemma~\ref{l:usep4}, we can, in any
case, use the pegs of $D \setminus P$ to get to a neighbor $t'$ of $t$
in at most $3$ moves.  Then, by Lemma~\ref{l:usep4}, we can use the pegs
of $P$ to get to a neighbor of $t'$ in at most $3$ more moves, after
which we can jump over $t'$ to $t$, thereby reaching $t$ in a total of
at most $7$ moves. Thus we may assume that all extended subgraphs coming
from $D \setminus P$ have independence number $2$, that the components
of $E'$ consist of the subgraphs $\check{C}$ for $C$ a component of $D
\setminus P$, and that no component $C$ of $D \setminus P$ has $|C| > 3$
and $|C| \ge 2\alpha(C) - 1$.

Therefore, the only components $C$ of $D \setminus P$ with $|C| -
2\alpha(\check{C}) \ge -1$ have order at most $3$ and, therefore, are
isomorphic to $K_1$, $K_2$, $P_3$, or $K_3$. Notice that $|C| -
2\alpha(\check{C})$ is $0$ for $C \cong K_2$ and is $-1$ for the other
possibilities for $C$ of order at most $3$; in particular, this
difference is never positive. Now, as $D \setminus P$ has at least
$2\alpha - 3$ vertices, a calculation analogous to that producing
Equation~(\ref{e:comps}) gives 
\[ 
\sum_{C \text{ a component of } D
\setminus P} (|C| - 2\alpha(\check{C})) \ge -1.  
\] 
Thus, for any subgraph $P$ of $D$ isomorphic to $P_4$, all components of
$D \setminus P$ are isomorphic to $K_2$, except for at most one
component isomorphic to $K_1$, $P_3$, or $K_3$. In particular, this
applies to all components of $D$ besides $K$. Therefore,
Equation~(\ref{e:comps}) yields 
\begin{equation}\label{e:k} 
|K| \ge 2\alpha(K) + 3.  
\end{equation} 
We now have three cases, based on the order of $K$.

\medskip\noindent\emph{Case 1:  We have $|K| \le 6$.}
Then $\alpha(K)$ is $1$, $|K|$ is $5$ or $6$, and $K$ is complete. Label
five of the vertices of $K$ as $v_1, \dots, v_5$. As before, some vertex
of $K$, say $v_1$, has distance $3$ from $t$, and we can take a path
$v_1 u w t$ of length $3$ from $v_1$ to $t$. Suppose that $u$ is not
adjacent to any vertex of $E \setminus K$. Then $K \cup \{u\}$ is still
a component of $E \cup \{u\}$, and $u$ is not adjacent to $t$, so we can
apply the argument leading to Equations~(\ref{e:comps}) and~(\ref{e:k}),
with $E$ replaced by $E \cup \{u\}$ and $\hat{K} = K \cup \{u\}$, to
obtain $|K| \ge 2\alpha(K \cup \{u\}) + 3$. This gives $\alpha(K \cup
\{u\}) = 1$, so that $K \cup \{u\}$ is complete. Thus we may reach $t$
in $4$ moves with the sequence $\peg{v_2}{v_1}{u}, \peg{v_3}{u}{w},
\peg{v_5}{v_4}{u}, \peg{u}{w}{t}$. Therefore, we may assume that $u$ is
adjacent to a vertex $x$ of $E \setminus K$.  Then we can move
$\peg{v_2}{v_1}{u}$ and $\peg{x}{u}{w}$ (the latter preceded by another
move to get a peg to $x$ if we had $x = u_C$), and then we get to $t$
via the sequence $\peg{v_4}{v_3}{v_1}, \peg{v_5}{v_1}{u},
\peg{u}{w}{t}$, reaching $t$ in at most $6$ moves, as desired.

\medskip\noindent\emph{Case 2:  We have $|K| = 7$.}
Then $\alpha(K)$ is at most $2$.  First suppose that $K$ is spanned by a
disjoint union of $P_2$ and $P_5$, with the endpoints of the $P_5$
non-adjacent.  Label the vertices of the $P_2$ as $v$ and $v'$ and the
vertices of the $P_5$ as $v_1, \dots, v_5$, with $v_i$ and $v_{i+1}$
adjacent, and $v_1$ and $v_5$ non-adjacent. Now either $v$ or $v'$, say
$v$, has distance $3$ from $t$, and we can take a path $v u w t$ of
length $3$ from $v$ to $t$. If $u$ is not adjacent to any vertex of $E
\setminus K$, then $K \cup \{u\}$ is a component of $E \cup \{u\}$, so,
as in Case~1, we get $|K| \ge 2\alpha(K \cup \{u\}) + 3$ and $\alpha(K
\cup \{u\}) = 2$, and $u$ must be adjacent to either $v_1$ or $v_5$, say
$v_5$. Therefore, in any event, $u$ must be adjacent to some vertex $x$
of $E \setminus \{v, v', v_1, \dots, v_4\}$. We can now move
$\peg{v'}{v}{u}$ and $\peg{x}{u}{w}$ (the latter preceded by another
move if we had $x = u_C$) and then use $v_1, \dots, v_4$ to reach $t$ in
$4$ more moves, for a total of at most $7$ moves. Thus we may assume
that $K$ is not spanned by a disjoint union of $P_2$ and $P_5$, with the
endpoints of the $P_5$ non-adjacent.

Let $P$ be a path in $K$ of maximal length.  Suppose that $P$ omitted at
least one vertex of $K$.  If the endpoints of $P$ were not adjacent,
then any vertex $v \in K \setminus P$ would be adjacent to one of them
(since we have $\alpha(K) \le 2$), and we could extend $P$. If the
endpoints of $P$ were adjacent, then the vertices of $P$ would form a
cycle, some vertex $v \in K \setminus P$ would be adjacent to a vertex
of $P$ (since $K$ is connected), and, again, we could extend $P$.  In
either case, we contradict the maximality of the length of $P$, so $P$
cannot omit any point of $P$ (and is, therefore, a Hamiltonian path of
$K$).

Label the vertices of $P$ (and, therefore, of $K$) as $w_1, \dots, w_7$,
with $w_i$ and $w_{i+1}$ adjacent.  Since $K$ is spanned by the union of
the two paths $w_1 w_2$ and $w_3 \dots w_7$, we must have $w_3$ and
$w_7$ adjacent, and, similarly, $w_1$ and $w_5$ must be adjacent. Next,
the paths $w_3 w_4$ and $w_2 w_1 w_5 w_6 w_7$ show that $w_2$ and $w_7$
are adjacent, as are $w_1$ and $w_6$. Now the paths $w_2 w_7$ and $w_3
w_4 w_5 w_1 w_6$ show that $w_3$ and $w_6$ are adjacent, as are $w_2$
and $w_5$. Finally, the paths $w_3 w_4$ and $w_1 w_2 w_5 w_6 w_7$ show
that $w_1$ and $w_7$ are adjacent.  Now one of $w_1$ and $w_7$, say
$w_1$, has distance $3$ from $t$, and we can take a path $w_1 u w t$ of
length $3$ from $w_1$ to $t$.  Finally, we can reach $t$ in $6$ moves
via the sequence $\peg{w_7}{w_1}{u}, \peg{w_3}{w_2}{w_1},
\peg{w_1}{u}{w}, \peg{w_4}{w_5}{w_1}, \peg{w_6}{w_1}{u}, \peg{u}{w}{t}$,
as desired.

\medskip\noindent\emph{Case 3:  We have $|K| \ge 8$.}
Let $P = v_1, \dots, v_4$ be a path on $4$ vertices in $K$.  We know
that all components of $K \setminus P$ are isomorphic to $K_2$, except
for at most one component isomorphic to $K_1$, $P_3$, or $K_3$.  Suppose
that one of $v_1$ and $v_2$ is adjacent to a vertex of a component of $K
\setminus P$ of order at least $2$, and that the same holds for one of
$v_3$ and $v_4$. There are at least two components of $K \setminus P$ of
order at least $2$ (since $K$ has order at least $8$), and each one is
adjacent to a vertex of $P$ (since $K$ is connected), so we can pick
distinct such components $C$ and $C'$ with $v_1$ or $v_2$ adjacent to a
vertex of $C$ and $v_3$ or $v_4$ adjacent to a vertex of $C'$. Then the
subgraph of $K$ induced by $\{v_1, v_2\} \cup C$ contains a subgraph
$P'$ isomorphic to $P_4$, and $K \setminus P'$ has a component of order
at least $4$ (namely, the one containing $\{v_3, v_4\} \cup C'$),
contradicting the assumption about the complement of any $P_4$ in $K$.
Thus, by symmetry, we may assume that neither $v_3$ nor $v_4$ is
adjacent to a vertex of a component of $K \setminus P$ of order at least
$2$.

Suppose that no component of $D \setminus P$ is isomorphic to $K_3$.
Then we have 
\begin{align*}
2\alpha(D) - 1 &  =  2\alpha(D) + 2\alpha(\{t\}) - 3\\
               & \le 2\alpha(G) - 3\\
               & \le |D| - 4\\
               &  =  |D \setminus P|\\
               &  =  \sum_{C \text{ a component of } D \setminus P} |C|\\
               & \le \sum_{C \text{ a component of } D \setminus P}
                          2\alpha(C)\\
               &  =  2\alpha(D \setminus P),
\end{align*}
giving $\alpha(D) = \alpha(D \setminus P)$.  This implies that every
vertex of $P$ is adjacent to a vertex of $D \setminus P$.  Applying this
to $v_3$ and $v_4$ shows that $D \setminus P$ must have a component
isomorphic to $K_1$, consisting of a single vertex $x$, and that both
$v_3$ and $v_4$ must be adjacent to this vertex $x$. Now the vertices
$v_1$ and $v_2$, together with the vertices of some component of $K
\setminus P$ isomorphic to $K_2$ (which must have a vertex adjacent to
either $v_1$ or $v_2$, since $K$ is connected), form a subgraph of $K$
isomorphic to $P_4$, whose complement in $K$ contains the complete
subgraph induced by $\{v_3, v_4, x\}$. Thus, replacing $P$ by this new
path if necessary, we may assume that $D \setminus P$ has a component
isomorphic to $K_3$.  In particular, no component of $D \setminus P$ is
isomorphic to either $K_1$ or $P_3$, and neither $v_3$ nor $v_4$ is
adjacent to any vertex of $K \setminus P$.

Let $\{x, y\}$ be any component of $K \setminus P$ isomorphic to $K_2$.
If $v_1$ were adjacent to, say, $x$, then the complement of the path $x
v_1 v_2 v_3$ in $K$ would contain the two isolated vertices $y$ and
$v_4$, contradicting the assumption about the complement of any $P_4$ in
$K$.  Therefore, $v_1$ is not adjacent to any vertex of a component of
$K \setminus P$ isomorphic to $K_2$.  Thus, $x$, say, is adjacent to
$v_2$.  The vertex $v_1$ cannot be adjacent to any vertex of a component
of $K \setminus P$, which would have to be the component isomorphic to
$K_3$, since otherwise, the complement of the path $y x v_2 v_3$ in $K$
would contain a component of order at least $4$ (that containing $v_1$
and the $K_3$). Finally, $v_2$ must be adjacent to a vertex of some
other component $C$ of $K \setminus P$ besides $\{x, y\}$, and the
complement in $K$ of a path on $x, v_2$, and two vertices of $C$
contains the singleton component $\{y\}$ and the component containing
$v_1$, which equals either $\{v_1\}$ or $\{v_1, v_3, v_4\}$.  This final
contradiction of the condition on complements of $P_4$ in $K$ shows
that, given our previous assumptions, Case 3 cannot occur, and this
proves the theorem.  
\end{proof}

\section{Conclusion}
We have introduced two new pegging quantities, namely the pegging number
and optimal pegging number of a graph.  We have successfully computed
these numbers for many classes of graphs, including paths, cycles,
joins, and products with complete graphs, using diverse tools including
basic pegging lemmas and pebbling.  In a forthcoming paper,
Wood~\cite{wood} studies pegging numbers of graph powers and products,
develops new general lower bounds for the pegging number, studies the
size of the reach of a distribution, and classifies some pegging moves
as unnecessary.

Any progress in the computation of these numbers and in the development
of computation tools would be of interest.  In particular, what is the
connection between the pegging numbers and other graph invariants such
as girth and connectivity?  Also, do pegging numbers behave nicely under
other graph operations like graph composition?  Are there any more
connections between pegging and pebbling? It should be noted that the
pegging analogue to Graham's conjecture~\cite{chung} is false, as shown
in~\cite{wood}.

In Theorem~\ref{t:reach}, we considered the effect of allowing stacking
and pebbling moves in pegging.  Let us refer to these two types of moves
as \emph{peggling} moves.  Let the peggling number (respectively, the
optimal peggling number) of a graph by the smallest positive integer $d$
such that every (respectively, some) multi-distribution $D$ of size $d$
has $\R_a(D) = V(G)$. 

Because peggling is like pebbling with more moves allowed, the peggling
and optimal peggling numbers of a graph are at most its pebbling and
optimal pebbling numbers, respectively. On the other hand,
Theorem~\ref{t:reach} shows that allowing stacking and pebbling moves
does not help reach additional targets in pegging.  Thus, because in
peggling the starting configuration can be any multi-distribution, the
peggling number of a graph is at least its pegging number, and the
optimal peggling number of a graph is at most its optimal pegging
number.  All of these inequalities may be strict; for example the graph
$P_4$ has pebbling number $8$, pegging number $3$, and peggling number
$5$; and the graph obtained by adding a pendant edge to each leaf of
$K_{1,3}$ has optimal pebbling number $4$, optimal pegging number $4$,
and optimal peggling number $3$.  Many of the basic properties of
pebbling, pegging, optimal pebbling, and optimal pegging numbers carry
over to peggling and optimal peggling numbers, and it would be
interesting to study these new quantities in more detail.

\section{Acknowledgments}

This research was done primarily at the University of Minnesota Duluth
Research Experience for Undergraduates summer progam with the support of
NSF grant DMS-9820438.  A special thanks is due to Joseph A. Gallian for
his encouragement and support, and we would also like to thank Melanie
Wood for comments on this paper.  Finally, we would like to thank the
participants of the 1994 University of Minnesota Duluth REU, who spent a
weekend in early investigations into pegging numbers and who introduced
the optimal pegging number.

\end{document}